

\input amstex
\documentstyle{amsppt}

\input label.def
\input degt.def

\input epsf
\def\picture#1{\epsffile{#1-bb.eps}}
\def\plot#1{\centerline{\epsfxsize.8\hsize\picture{plot-#1}}}

\def\ie{\emph{i.e.}}
\def\eg{\emph{e.g.}}
\def\cf.{\emph{cf}\.}
\def\via{\emph{via}}
\def\etc{\emph{etc}}

{\catcode`\@11
\gdef\proclaimfont@{\sl}}

\def\dash{\item"\hfill--\hfill"}
\def\Dashes{\widestnumber\item{--}}

\loadbold
\def\bA{\bold A}
\def\bD{\bold D}
\def\bE{\bold E}
\def\bJ{\bold J}

\let\splus\oplus

\let\Ga\alpha
\let\Gb\beta
\let\Gg\gamma
\let\Gd\delta
\let\Ge\epsilon
\let\Gs\sigma
\let\Gr\rho
\let\Gf\varphi
\def\1{^{-1}}
\def\ba{\bar a}
\def\bc{\bar c}
\def\bGa{\bar\Ga}
\def\bGg{\bar\Gg}

\def\bSigma{\smash{\bar\Sigma}}

\def\BB{\bar B}
\def\PP{\bar P}
\def\RR{\bar R}
\def\LL{\bar L}
\def\OO{\bar O}

\def\tG{\tilde G}
\def\tA{\tilde A}
\def\tX{\tilde X}
\def\tSigma{\tilde\Sigma}
\def\tc{\tilde c}
\let\L\Lambda

\def\bp{\bar p}
\def\bq{\bar q}
\def\CO{\Cal O}

\def\CG#1{\Z_{#1}}     
\def\BG#1{\Bbb B_{#1}} 
\def\SG#1{\Bbb S_{#1}} 
\def\AG#1{\Bbb A_{#1}} 
\def\DG#1{\Bbb D_{#1}} 

\def\F{\Bbb F}
\def\Cp#1{\Bbb P^{#1}}
\def\CK{\Cal K}
\def\CL{\Cal L}

\def\Sk{\operatorname{Sk}}
\def\discr{\operatorname{discr}}
\def\Aut{\operatorname{Aut}}

\let\<\langle
\let\>\rangle
\def\ls|#1|{\mathopen|#1\mathclose|}

{\let\1\gdef\let\+\relax
\gdef\setcat{\catcode`\a\active \catcode`\d\active \catcode`\e\active
\catcode`\+\active}
\setcat
\1a[#1]{\bA_{#1}} \1d[#1]{\bD_{#1}} \1e[#1]{\bE_{#1}}
}
{\catcode`\+\active\gdef+{\splus}}
\def\getline#1,#2[#3],#4[#5],#6,#7[#8]{#2&#1\cr}
{\obeylines%
\gdef\TAB#1{\vtop\bgroup\openup1pt\obeylines\let
\getline\setcat\expandafter\halign\expandafter\bgroup#1\cr}}
\def\ENDTAB{\crcr\egroup\egroup}
\def\Torus{\llap{$##$\ }&$##$\hss}
\def\NonTorus{$##$\hss&$##$\hss}

\topmatter

\author
Alex Degtyarev
\endauthor

\title
Irreducible plane sextics\\with large fundamental groups
\endtitle

\address
Department of Mathematics,
Bilkent University,
06800 Ankara, Turkey
\endaddress

\email
degt\@fen.bilkent.edu.tr
\endemail

\abstract
We compute the fundamental group of the complement of each
irreducible sextic of weight eight or nine (in a sense, the
largest groups for irreducible sextics), as well as of $169$ of
their derivatives (both of and not of torus type). We also give a
detailed geometric description of sextics of weight eight and nine
and of their moduli spaces and compute their Alexander modules;
the latter are shown to be free over an appropriate ring.
\endabstract



\endtopmatter

\document

\section{Introduction\label{S.intro}}

\subsection{Sextics of torus type}
Recall that a plain sextic~$B$ is said to be of \emph{torus type}
(more precisely, \emph{$(2,3)$-torus type}) if its equation can be
represented in the form $p^3+q^2=0$, where $p$ and~$q$ are some
homogeneous polynomials of degree~$2$ and~$3$, respectively.
A singular point~$P$ of~$B$ is called \emph{inner} (with respect
to a given torus structure $(p,q)$\,) if it belongs to the conic
$\{p=0\}$; then $P$ also belongs to the cubic $\{q=0\}$.
Otherwise, $P$ is called \emph{outer}.
A sextic
of torus type is called \emph{tame} if all its singular points are
inner.

In spite of the algebraic definition, the property of being of
torus type is of a purely topological nature; in particular, it is
invariant under equisingular deformations. For example, an
irreducible sextic is of torus type if and only if its Alexander
polynomial (see A.~Libgober~\cite{Libgober1} and
Definition~\ref{def.Alexander} below) is nontrivial,
see~\cite{degt.Oka} and~\cite{degt.Oka2}; the latter condition
can be restated as
$\ls|\pi'\!/\pi''|=\infty$, where $\pi=\pi_1(\Cp2\sminus B)$ and
$'$ stands for the derived group.
Alternatively, an irreducible
sextic is of torus type if and only if its fundamental group
factors to the dihedral group~$\DG6$,
see~\cite{degt.Oka} and~\cite{degt.Oka2} or
H.~Tokunaga~\cite{Tokunaga.new} (where reducible sextics with
simple singularities are considered as well).

Define the \emph{weight} $w(P)$ of a simple singular point~$P$ as
follows: $w(\bA_{3k-1})=k$, $k\in\Z$, $w(\bE_6)=2$, and $w(P)=0$
for all other types. The \emph{weight} $w(B)$ of a sextic~$B$ with
simple singularities is the total weight of its singular
points. According to~\cite{degt.Oka}, any sextic of torus type has
weight $6\le w\le9$, and any irreducible sextic of weight $w\ge7$
is of torus type. All sets of singularities of maximal weight~$8$
or~$9$ are classified in~\cite{degt.Oka}: each sextic of weight
nine has nine cusps (and is dual to a nonsingular cubic),
and sextics of weight eight have one of the
following sets of singularities:
$$
\gathered
\bE_6\splus\bA_5\splus4\bA_2,\quad
\bE_6\splus6\bA_2,\quad
2\bA_5\splus4\bA_2,\\
\bA_5\splus6\bA_2\splus\bA_1,\quad
\bA_5\splus6\bA_2,\quad
8\bA_2\splus\bA_1,\quad
8\bA_2.
\endgathered\eqtag\label{list.w=8}
$$
In addition, it seems reasonable to assign weight eight to the
following two non-simple sets of singularities
(see
the discussion on the Alexander polynomial below):
$$
\bJ_{2,3}\splus3\bA_2,\quad
\bJ_{2,0}\splus4\bA_2.\eqtag\label{list.J}
$$
(We use Arnol$'$d's notation~\cite{Arnold} for non-simple
singularity types: $\bJ_{2,p}$ is a point of simplest
tangency of a smooth
branch and a singularity of type~$\bA_{p+3}$, $p\ge0$.) The
lists~\eqref{list.w=8} and~\eqref{list.J} appeared first in
M.~Oka~\cite{Oka.poly}: they are the sets of singularities realizing the
Alexander polynomial $(t^2-t+1)^2$.

From the results of Oka~\cite{Oka.poly}, it follows that
the Alexander polynomial $\Delta_B(t)$ of an
irreducible sextic~$B$ of torus type is $(t^2-t+1)^s$, where $s=1$
if $w(B)=6$ and $s=w(B)-6$ if $w(B)\ge7$. (A simple proof of the
fact that $s=1$ for $w(B)\le7$ is given in
Section~\ref{s.perturbations}, see Corollary~\ref{cor.Alexander}.)
Among irreducible sextics with a non-simple singular point, only
the two listed in~\eqref{list.J} have Alexander polynomial
$(t^2-t+1)^2$ (see~\cite{Oka.poly} or~\cite{degt.Oka}
and~\cite{degt.Oka2}); that is why they are assigned weight
eight.

\subsection{Principal results}
Since $\deg\Delta_B=\dim_\C(\pi'\!/\pi'')\otimes\C$, where
$\pi=\pi_1(\Cp2\sminus C)$, one can conclude that, in a sense, the
groups of sextics of weight eight and nine are largest
possible for irreducible sextics. The principal result of the
present paper is the computation of these groups.

\theorem\label{th.w=8}
The fundamental group $\pi_1(\Cp2\sminus B)$ of a plane sextic~$B$
with a set of singularities~$\Sigma$ as in~\eqref{list.w=8}
or~\eqref{list.J} is as follows\rom:
\Dashes
\roster
\dash
the group~$G_2$ given by~\eqref{eq.G2}, if
$\Sigma=\bE_6\oplus\bA_5\oplus4\bA_2$\rom;
\dash
the group~$G_1$ given by~\eqref{eq.G1}, if
$\Sigma=\bA_5\oplus6\bA_2\oplus\bA_1$\rom;
\dash
the group~$G_0$ given by~\eqref{eq.G0} otherwise.
\endroster
\endtheorem

The curves with the set of singularities $\bJ_{2,3}\splus3\bA_2$
are tame; their fundamental group~$G_0$ is found in M.~Oka,
D.~T.~Pho~\cite{OkaPho}.

\remark{Remark}
There are two obvious
perturbation epimorphisms
$G_2\twoheadrightarrow G_1\twoheadrightarrow G_0$. The latter is
proper, see Corollary~\ref{no.iso}. At present, I do not know
whether the former epimorphism $G_2\twoheadrightarrow G_1$
is proper or not.
\endremark

\theorem\label{th.w=9}
The fundamental group $\pi_1(\Cp2\sminus B)$ of a plane sextic~$B$
with the set of singularities~$9\bA_2$ is the group~$G_3$ given
by~\eqref{eq.G3}.
\endtheorem

Theorems~\ref{th.w=8} and~\ref{th.w=9} are proved in
Sections~\ref{s.w=8} and~\ref{s.w=9}, respectively. The
fundamental group of the nine-cuspidal sextic
(Theorem~\ref{th.w=9}) was first computed by
Zariski~\cite{Zariski.9a2}; yet another computation and further
generalizations can be found in
J.~I.~Cogolludo~\cite{Cogolludo}.

As a necessary preliminary step, we describe an explicit geometric
construction of the curves in question and their torus structures,
see Section~\ref{s.moduli}, and study their moduli space. (The key
ingredient here is Theorem~\ref{th.reduction}, stating that each
sextic of weight eight or nine is symmetric.)
In particular,
we prove the uniqueness of an equisingular deformation family
realizing each of the sets of singularities listed
in~\eqref{list.w=8}.
(For~\eqref{list.J} and for $9\bA_2$
the uniqueness is known.)

\theorem\label{th.1}
Each set of singularities listed in~\eqref{list.w=8} is realized by
a single equisingular deformation family of plane sextics.
\endtheorem

Theorem~\ref{th.1} admits a refinement, taking into account the
torus structure.

\theorem\label{th.torus}
With one exception, a pair $(B,\{\text{\rm torus structure}\})$,
where~$B$ is an irreducible
sextic of weight eight,
is determined up to equisingular deformation of torus
structures by the combinatorial type of the pair
$(\Sigma,\Sigma_{\inj})$,
where $\Sigma$ is the set of singularities of~$B$ and
$\Sigma_{\inj}$ is the set of its
inner singularities.
The exception is the pair
$(\Sigma,\Sigma_{\inj})=(\bE_6\splus\bA_5\splus4\bA_2,\bE_6\splus\bA_5\splus2\bA_2)$,
which is realized by two complex conjugate
equisingular deformation families.
\endtheorem

It should not be very difficult to deduce Theorems~\ref{th.1}
and~\ref{th.torus} arithmetically, using~\cite{JAG}. However, we
give a geometric proof, see Section~\ref{proof.1}, based on a
detailed description of the curves and their moduli
space.

In Section~\ref{s.torus.9a2}, we discuss geometric
properties of the twelve torus structures of a sextic of weight
nine and prove an analogue of Theorem~\ref{th.torus}.

As a first application of Theorems~\ref{th.w=8} and~\ref{th.w=9},
we compute the Alexander modules (in the
sense of Libgober~\cite{Libgober2}, see
Definition~\ref{def.Alexander} below) of all sextics of weight
eight and nine. The following statement is proved
in Section~\ref{proof.Alexander}.

\theorem\label{th.Alexander}
The Alexander module of a sextic of weight eight
\rom(nine\rom) is a free module on two \rom(respectively,
three\rom) generators over the ring $\L=\Z[t]/(t^2-t+1)$.
\endtheorem

As another application, we compute the fundamental groups of the
perturbations of sextics of weight eight. Altogether,
we consider $47$ sets of singularities of torus
type, see Theorem~\ref{th.w=6,7},
and $122$ sets of singularities that are not of torus type
and not covered by M.~V.~Nori's theorem~\cite{Nori}, see
Theorem~\ref{th.nontorus}. For about half of these sets of
singularities, the uniqueness of an equisingular deformation
family is known. Among them are $17$ so called \emph{classical
Zariski pairs}, see Corollary~\ref{cor.Zariski}; we show that, as
in the original example by O.~Zariski~\cite{Zariski.group},
the fundamental groups of the two curves constituting
each pair are $\CG2*\CG3$ and~$\CG6$.

It is worth mentioning that irreducible sextics of torus type,
their moduli, and fundamental groups have
been a subject of intensive study, so that some of the results of
this paper may overlap with results obtained by other authors. We
will cite Oka, Pho~\cite{OkaPho.moduli} (classification of
singularities of torus type
and moduli of maximal sextics of torus type), Oka,
Pho~\cite{OkaPho} (fundamental groups of tame sextics),
Oka~\cite{Oka.pairs} (Zariski pairs involving sextics of torus
types),
and, for
further references, recent paper Eyral, Oka~\cite{EyralOka}.

\subsection{Digression: geometry of a nine cuspidal sextic}
As a by-product, we apply the results obtained in the paper
to study the
geometry of a nine cuspidal sextic (sextic of weight nine). Any
such sextic~$B$ is an elliptic curve. We show that,
if $B$ is generic (no
complex multiplication), then the group $\Aut(\Cp2,B)$ of
projective automorphisms of~$B$ is a semi-direct product
$(\CG3\times\CG3)\rtimes\CG2$, see Section~\ref{s.Aut.9a2}. Then,
in Section~\ref{s.torus.9a2}, we characterize the twelve torus
structures of~$B$, see~\cite{degt.Oka} or~\cite{Tokunaga},
in terms of their
stabilizers, which are dihedral subgroups
$\DG6\subset\Aut(\Cp2,B)$, as well as in terms of the inflection
points of the dual cubic curve. Finally, we prove an analogue of
Theorem~\ref{th.torus}: pairs consisting of a sextic~$B$ of weight
nine and a torus structure on~$B$ form a connected deformation
family.

Some of these results may be known.

\subsection{Contents of the paper}
The starting point of our computation is~\S\ref{S.construction}:
we apply the characterization found in~\cite{degt.Oka},
V.~V.~Nikulin's theory of discriminant forms, and the theory of
periods of $K3$-surfaces to show that each sextic~$B$
of weight eight
or nine with simple singularities is symmetric
(Theorem~\ref{th.reduction}),
thus folding~$B$
to a very special trigonal curve~$\BB$ in Hirzebruch
surface~$\Sigma_2$ (quadratic cone). It appears that
the symmetry constructed is the only non-trivial projective
automorphism of a generic curve of weight eight. In the case of
weight nine, the situation is different:
we show that the automorphism group of a generic nine cuspidal
sextic has order~$18$.

In~\S\ref{S.BB}, we represent the trigonal curve~$\BB$ obtained
above by explicit equations and use them to analyze the
automorphisms of~$\BB$, its real structures, torus structures, and
special (tangent, double tangent, inflection tangent, \etc.)
sections. The results are applied to prove Theorems~\ref{th.1}
and~\ref{th.torus}, as well as to the study of generic
nine cuspidal sextics.

In~\S\ref{S.group}, we use the equations developed in~\S\ref{S.BB}
to visualize the braid monodromy of~$\BB$ and to write down
presentations for the fundamental groups. Theorems~\ref{th.w=8},
\ref{th.w=9}, and~\ref{th.Alexander} are proved here.

Finally, in~\S\ref{S.appl}, we use the presentations obtained above
to compute the groups of a number of other sextics.
The key ingredient here is Proposition~\ref{subgraph}, which
states that any `thinkable' perturbation of singularities of a
plane sextic is indeed realizable. We discuss briefly the
uniqueness of equisingular deformation families realizing the sets
of singularities for which the groups are found.

\section{The construction}\label{S.construction}

In this section, we show that any irreducible sextic of weight
eight or nine can be obtained as a double covering of a certain
trigonal curve~$\BB$ in a geometrically ruled rational
surface~$\Sigma_2$.

\subsection{Statements}
Recall that any involutive automorphism $c\:\Cp2\to\Cp2$ has a
fixed line~$L_c$ and an isolated fixed point~$O_c$, and the
quotient $\Cp2(O_c)\!/c$ is the rational geometrically ruled
surface~$\Sigma_2$.
The images in~$\Sigma_2$ of~$O_c$ and~$L_c$ are,
respectively, the exceptional section~$E$ and a generic
section~$\LL$, so that $\Cp2(O_c)$ is the double covering
of~$\Sigma_2$ ramified at $\LL+E$. Alternatively, $\Cp2$ is the
double covering of the quadratic cone $\Sigma_2/E$ ramified
at~$\LL$ and $E/E$. (Here and below, we use the notation
$\Cp2\!/c$ for the orbit space of~$c$, and $\cdot/E$ stands for
the contraction of~$E$, so that $E/E$ is the vertex of the cone
$\Sigma_2/E$.)

The principal result of this section is the following theorem.

\theorem\label{th.reduction}
Let $B\subset\Cp2$ be a sextic of weight eight or nine and with
simple singularities only.
Then $\Cp2$ admits an
involution~$c$ preserving~$B$, so that $O_c\notin B$,
and the image of~$B$ in $\Sigma_2=\Cp2(O_c)\!/c$
is a trigonal curve~$\BB$, disjoint from~$E$,
with four cusps~$\bA_2$. Conversely, given a trigonal curve
$\BB\subset\Sigma_2$ as above and a section~$\LL$ of~$\Sigma_2$
disjoint from~$E$, the double covering of~$\BB$ ramified at~$\LL$
and $E/E$ is a plane sextic of weight eight or nine and with
simple singularities only.
\endtheorem

Theorem~\ref{th.reduction} is proved at the end of this section,
in~\ref{s.proof} below.

\remark{Remark}
According to~\cite{degt.Oka2}, the two deformation families of
weight~$8$ with non-simple singular points, see~\eqref{list.J},
can also be obtained from a four cuspidal trigonal curve, but by a
birational transformation rather than double covering.
\endremark

\subsection{Discriminant forms}\label{s.discr}
An \emph{\rom(integral\rom) lattice} is a finitely generated free
abelian group~$L$ supplied with a symmetric bilinear form $b\colon
L\otimes L\to\Z$. We abbreviate $b(x,y)=x\cdot y$ and
$b(x,x)=x^2$. A lattice~$L$ is called \emph{even} if $x^2=0\bmod2$ for
all $x\in L$. As the transition matrix between two integral bases
has determinant $\pm1$, the determinant $\det L\in\Z$ (\ie, the
determinant of the Gram matrix of~$b$ in any basis of~$L$) is well
defined. A lattice~$L$ is called \emph{nondegenerate} if $\det
L\ne0$; it is called \emph{unimodular} if $\det L=\pm1$.

Given a lattice~$L$,
the bilinear form extends to a form
$(L\otimes\Q)\otimes(L\otimes\Q)\to\Q$.
If
$L$ is nondegenerate, the dual group $L^*=\Hom(L,\Z)$ can
be identified with the subgroup
$$
\bigl\{x\in L\otimes\Q\bigm|
 \text{$x\cdot y\in\Z$ for all $x\in L$}\bigr\}.
$$
In particular, $L\subset L^*$
is a finite index subgroup. The quotient $L^*\!/L$
is called the \emph{discriminant group} of~$L$ and is
denoted by $\discr L$ or~$\CL$. The discriminant group
inherits from $L\otimes\Q$ a symmetric bilinear form
$b_{\CL}\:\CL\otimes\CL\to\Q/\Z$,
called the \emph{discriminant form}, and,
if $L$ is even, its quadratic extension $q_{\CL}\:\CL\to\Q/2\Z$.
One has $\#\CL=\ls|\det L|$; in particular,
$\CL=0$ if and only if $L$ is unimodular.

From now on, \emph{all lattices considered are even}.

An \emph{extension} of a lattice~$L$ is another lattice~$M$
containing~$L$, so that the form on~$L$ is the restriction of that
on~$M$. An \emph{isomorphism} between two extensions
$M_1\supset L$ and $M_2\supset L$ is an isometry $M_1\to M_2$
whose restriction to~$L$ is the identity.
In what follows,
we are only interested in the case when
$[M:L]<\infty$. Next
two theorems are found in Nikulin~\cite{Nikulin}.

\theorem\label{th.Nik1}
Given a nondegenerate
lattice~$L$, there is a canonical one-to-one correspondence
between the set of isomorphism classes of finite index extensions
$M\supset L$
and the set of isotropic subgroups $\CK\subset\CL$.
Under this correspondence,
$M=\{x\in L^*\,|\,x\bmod L\in\CK\}$ and
$\discr M=\CK^\perp\!/\CK$.
\qed
\endtheorem

The isotropic subgroup $\CK\subset\CL$ as in Theorem~\ref{th.Nik1}
is called the \emph{kernel} of the extension $M\supset L$. It can
be defined as the image of $M\!/L$ under the homomorphism induced by
the natural inclusion $M\hookrightarrow L^*$.

\theorem\label{th.Nik2}
Let $M\supset L$ be a finite index
extension of a nondegenerate lattice~$L$,
and let $\CK\subset\CL$ be its
kernel.
Then, an auto-isometry $L\to L$ extends to~$M$ if and only if the
induced automorphism of~$\CL$ preserves~$\CK$.
\qed
\endtheorem

We will use Theorem~\ref{th.Nik2} in the following form.

\corollary\label{extension}
Let $L\subset M$ be a nondegenerate sublattice of a unimodular
lattice~$M$, and let $\CK\subset\CL$ be the kernel of the
extension $\tilde L\supset L$, where $\tilde L$ is the primitive
hull of~$L$ in~$M$. Consider an auto-isometry $c\:L\to L$. Then,
$c\oplus\id_{L^\perp}$ extends to~$M$ if and only if $c$
preserves~$\CK$ and the auto-isometry of
$\CK^\perp\!/\CK$ induced by~$c$ is the identity.
\endcorollary

\proof
We apply Theorem~\ref{th.Nik2} twice: first, to the extension
$\tilde L\supset L$, then to the extension
$M\supset\tilde L\oplus L^\perp$.
The condition that $c$ should preserve~$\CK$ is necessary and
sufficient for~$c$ to extend to an isometry~$\tilde c$
of~$\tilde L$. If it does extend, the auto-isometry of
the discriminant
$$
\discr(\tilde L\oplus L^\perp)=(\CK/\CK^\perp)\oplus\discr L^\perp
$$
induced by~$\tilde c\oplus\id_{L^\perp}$ is
$\tilde c_{\CK}\oplus\id_{\discr L^\perp}$, where $\tilde c_{\CK}$
is the automorphism induced by~$\tilde c$ (or~$c$) on
$\CK/\CK^\perp$. Since the kernel of the extension
$M\supset\tilde L\oplus L^\perp$ is the graph of a certain
anti-isometry $\CK/\CK^\perp\to\discr L^\perp$,
see Nikulin~\cite{Nikulin}, it is preserved by the automorphism
$\tilde c_{\CK}\oplus\id_{\discr L^\perp}$ above (the condition
necessary and sufficient for the auto-isometry
$\tilde c\oplus\id_{L^\perp}$ to
extend further to~$M$) if and only if
$\tilde c_{\CK}=\id_{\CK/\CK^\perp}$.
\endproof

\subsection{Sextics of weight eight}\label{s.c.w=8}
Let $B$ be a sextic with one of the sets of singularities listed
in~\eqref{list.w=8}. Denote by~$\tX$ the resolution of
singularities of the double covering $X\to\Cp2$ ramified at~$B$.
It is a $K3$-surface. For each singular point~$P$ of~$B$, let
$\Gamma_P$ be the incidence graph of the exceptional divisors
in~$\tX$ over~$P$
(it is a Dynkin diagram of the same name as the type of~$P$),
and let $\Sigma_P\subset H_2(\tX)$ be the
sublattice spanned by the vertices of~$\Gamma_P$. Let
$\Gamma=\bigcup_P\Gamma_P$ and $\Sigma=\bigoplus_P\Sigma_P$, and
denote by $\tSigma$ the primitive hull of~$\Sigma$ in $H_2(\tX)$.

For future references, introduce also the class $h\in H_2(\tX)$
realized by the pull-back of a generic line in~$\Cp2$. Observe
that, regarded as an element of the Picard group, $h$ is the
linear system defining the projection $\tX\to\Cp2$.

Split all singular points of~$B$ of positive weight
into four groups of total weight two each, and let $\Gamma_i$ and
$\Sigma_i$, $i=1,\ldots,4$, be the corresponding subgraphs
of~$\Gamma$ and sublattices of~$\Sigma$, respectively.
Each~$\Sigma_i$ is either $2\bA_2$ or~$\bA_5$ or~$\bE_6$. Let
$c_i\:\Sigma_i\to\Sigma_i$ be the automorphism induced by a
non-trivial symmetry of~$\Gamma_i$: in the cases~$\bA_5$
and~$\bE_6$, such a symmetry is unique; in the case~$2\bA_2$, there
are two symmetries transposing the two components, and we pick one
of them. There is a unique pair of $c_i$-skew-invariant
nonzero elements $\pm x_i\in\discr\Sigma_i$; one has
$x_i^2=\frac23\bmod2\Z$. (In the cases~$\bA_5$ and~$\bE_6$,
$\pm x_i$ are the generators of
$\discr\Sigma_i\otimes\F_3=\F_3$; in the
case~$2\bA_2$, there are two pairs of opposite elements of square
$\frac23\bmod2\Z$, and the choice of one of them
determines~$c_i$.) Choose one of them and
denote it by~$x_i$.

Now, the description of sextics of weight~$8$ given
in~\cite{degt.Oka} can be restated as follows: the
sublattices~$\Sigma_i$, involutions~$c_i$, and elements~$x_i$,
$i=1,\ldots,4$, as
above can be chosen so that the kernel $\CK\subset\discr\Sigma$ of
the extension $\tSigma\supset\Sigma$ is spanned by $x_1+x_2+x_3$
and $x_1-x_2+x_4$. Then, extending $\bigoplus_ic_i$ identically to
the rest of~$\Sigma$, we obtain an involution
$c_\Sigma\:\Sigma\to\Sigma$ with the following properties:
$c_\Sigma$ acts identically on the $p$-primary part of
$\discr\Sigma$ for any prime $p\ne3$, and $\CK$ is a maximal (half
dimension) isotropic subspace of the $(-1)$-eigenspace
of~$c_\Sigma$ in the $3$-primary part $\discr\Sigma\otimes\F_3$.
Hence, $\CK^\perp\!/\CK$ can be identified with the
$c_\Sigma$-invariant part of $\discr\Sigma$, and, due to
Corollary~\ref{extension},
$c_\Sigma\oplus\id_{\Sigma^\perp}$
extends to an involution $\tc_*$ on $H_2(\tX)$.

\subsection{Sextics of weight nine}\label{s.c.w=9}
Let~$B$ be a sextic with the set of singularities~$9\bA_2$. Pick
one of the cusps and treat it as an ordinary point (of weight
zero). Applying the construction of Section~\ref{s.c.w=8}, we
obtain a splitting of the remaining eight cusps into four groups
and an involution $c_\Sigma\:\Sigma\to\Sigma$. According
to~\cite{degt.Oka}, the kernel~$\CK$ is spanned by the elements
$x_1+x_2+x_3$ and $x_1-x_2+x_4$ introduced above and by
$y_1+y_2+y_3+y_4+z_9$, where $y_i$ is an appropriately chosen
generator of~$x_i^\perp$ in $\discr\Sigma_i$, $i=1,\ldots,4$, and
$z_9$ is a generator of the discriminant of the ninth cusp.
Thus, $\CK^\perp\!/\CK$ can be identified with
a certain subquotient of the the
$c_\Sigma$-invariant part of $\discr\Sigma$, and still
$c_\Sigma\oplus\id_{\Sigma^\perp}$
extends to an involution $\tc_*$ on $H_2(\tX)$.

\subsection{Proof of Theorem~\ref{th.reduction}}\label{s.proof}
Consider the $K3$-surface~$\tX$ introduced
in~\ref{s.c.w=8} and the involution~$\tc_*$ on $H_2(\tX)$
constructed in~\ref{s.c.w=8} and~\ref{s.c.w=9}. From the
construction, it follows that $\tc_*$ preserves~$\Gamma$ (as a
set), the class $h\in H_2(\tX)$ of
the pull-back of a line, see~\ref{s.c.w=8},
and the
class $\omega\in H_2(\tX;\C)$ of a holomorphic $2$-form on~$\tX$ (as
both~$h$ and~$\omega$ are orthogonal to~$\Sigma$). Then, $\tc_*$
also preserves the positive cone of~$\tX$ (as it can be expressed
in terms of~$\Gamma$, $h$, and~$\omega$, \cf.~\cite{degt.Oka3})
and hence $\tc_*$ is induced by a unique involution
$\tc\:\tX\to\tX$; the latter is symplectic (\ie, preserving
holomorphic $2$-forms) and commutes with the deck translation of
the covering $\tX\to\Cp2$ (as $\tc_*$ preserves~$h$). The descent
of~$\tc$ to~$\Cp2$ is the desired involution $c\:\Cp2\to\Cp2$.

The image~$\BB$ of~$B$ in $\Sigma_2=\Cp2(O_c)\!/c$ can easily be
studied similar to~\cite{degt.Oka3}. (In particular, all necessary
facts relating the singularities of~$B$ and those of $\BB+\LL$ are
found in~\cite{degt.Oka3}; the relevant part is represented in
Table~\ref{tab.data} in \S\ref{S.BB}.)
For the converse statement, one
observes that each cusp~$\PP$ of~$\BB$ gives rise to either two
cusps of~$B$ (if $\PP\notin\LL$) or one singular point of~$B$ of
type~$\bA_5$ or~$\bE_6$ (if $\PP\in\LL$ and, for~$\bE_6$, $\LL$ is
tangent to~$\BB$ at~$\PP$). The intersection points $\LL\cap\BB$
nonsingular for~$\BB$ produce type~$\bA$ singularities of~$B$.
Hence, all singularities of~$B$ are simple and the total weight
of~$B$ is at least eight.
\qed

\subsection{Automorphisms of nine cuspidal sextics}\label{s.Aut.9a2}
Let $B\subset\Cp2$ be a nine cuspidal sextic. As follows
from~\ref{s.c.w=9} and~\ref{s.proof}, for each cusp~$P$ of~$B$,
there is an involution of the pair $(\Cp2,P)$ preserving~$P$ and
transposing the other eight cusps. This classical fact has a
transparent geometric explanation. Consider the nonsingular cubic
$C\subset\Cp2$ dual to~$B$. Then, the automorphisms of $(\Cp2,B)$
are those of $(\Cp2,C)$, which in turn are the automorphisms of
the abstract elliptic curve~$C$ preserving its nine inflection
points. (Note that $C$ is the normalization of~$B$;
as abstract curves they can be identified.)
Hence, the following statement holds.

\lemma\label{9a2.Aut}
If $B\cong C$ is generic \rom(no complex multiplication\rom), then
the
group $\Aut(\Cp2,B)$ of projective automorphisms
of~$B$
is a semi-direct product of $\CG3\times\CG3$
\rom(affine shifts by the vectors
$(\frac m3,\frac n3)\in B\cong\R/\Z\times\R/\Z$\rom) and~$\CG2$
\rom(multiplication by~$(-1)$ in~$B$\rom),
$\CG2$ acting on $\CG3\times\CG3$
\via\ $z\mapsto-z$.
\qed
\endlemma

The nine involutions mentioned above are the nine
elements of the non-trivial coset modulo
$\CG3\times\CG3$. These are all order~$2$ elements in
$\Aut(\Cp2,B)$.

The automorphisms of~$C$ can be seen from its Hesse's normal form
$$
z_0^3+z_1^3+z_2^3=\lambda z_0z_1z_2.
$$
The group $\Aut(\Cp2,C)$ is generated by the permutations of the
coordinates and the automorphisms
$(z_0:z_1:z_2)\mapsto(z_0:\Ge z_1:\Ge^2z_2)$, $\Ge^3=1$.

If $C$ has a complex multiplication, then $\Aut(\Cp2,C)$ is a
semi-direct product of $\CG3\times\CG3$ and~$\CG4$ or~$\CG6$; in
particular, the set of involutions in $\Aut(\Cp2,C)$ is the same.
As a consequence, we have the following statement.

\lemma\label{9a2.involutions}
For a nine cuspidal sextic $B\subset\Cp2$, each involutive
automorphism $c\in\Aut(\Cp2,B)$ has a single fixed cusp~$P_c$, and
$c$ is determined by~$P_c$ uniquely.
\qed
\endlemma

The action of $\Aut(\Cp2,B)$ on the torus structures of~$B$
is discussed in~\ref{s.torus.9a2} below.

\section{The trigonal curve $\BB$}\label{S.BB}

In this section, we study the geometry of the four cuspidal
trigonal curve $\BB\subset\Sigma_2$. Most calculations involving
equations were performed using {\tt Maple}.

\subsection{The geometric description}\label{s.BB}
A trigonal curve $\BB\subset\Sigma_2$ with
the set of singularities $4\bA_2$
is a maximal trigonal curve in the sense
of~\cite{degt.kplets}; in particular, it is unique up to
automorphism of~$\Sigma_2$. The skeleton $\Sk\subset\Cp1$ of~$\BB$
(see~\cite{degt.kplets}) is the $1$-skeleton of a regular
tetrahedron~$\Delta$ (assuming that $\Cp1$ is regarded as the
surface of~$\Delta$). Hence, the group of Klein automorphisms
of~$\BB$ (\ie, holomorphic or anti-holomorphic automorphisms
$\Sigma_2\to\Sigma_2$ preserving~$\BB$) is the full symmetric
group~$\SG4$ (the group of symmetries of~$\Delta$), and the group
of holomorphic automorphisms of~$\BB$ is the subgroup
$\AG4\subset\SG4$ (the group of rotations of~$\Delta$). Both
groups act faithfully on the set of cusps of~$\BB$ (barycenters of
the faces of~$\Delta$). Explicit generators for the group~$\AG4$
are given in~\eqref{eq.order3} and~\eqref{eq.order2} below.

As a consequence, $\BB$ is real with respect to six different real
structures on~$\Sigma_2$ (transpositions in~$\SG4$).
Each real structure preserves exactly
two of the four cusps of~$\BB$ and is uniquely determined by this
pair of cusps. Below, we use two distinct real structures,
see~\eqref{eq.conj}, to
visualize the braid monodromy of~$\BB$.

Alternatively, $\BB$ can be constructed as a birational transform
of a three cuspidal plane quartic~$C$. As is known, $C$ has a
unique double tangent~$L$; one should pick one of the tangency
points~$P$, blow it up twice, and blow down the proper transform
of~$L$ and one of the exceptional divisors over~$P$. The inverse
transformation is the stereographic projection of the quadratic
cone $\Sigma_2/E$ from one of the cusps of~$\BB$.

\subsection{The equation of~$\BB$}\label{s.equation}
In appropriate affine coordinates $(x,y)$ in~$\Sigma_2$ a
curve~$\BB$ as above can be given by the polynomial
$$
f(x,y)=4y^3-(24x^3+3)y+(8x^6+20x^3-1).\eqtag\label{eq.equation}
$$
The discriminant of this expression with respect to~$y$ is
$$
-(2)^{10}(3)^3x^3(x^3-1)^3;\eqtag\label{eq.discrim}
$$
hence, $\BB$ does have four cusps.

The curve is rational; it can be parametrized by
$$
x(t)=\frac{3t}{t^3+2},\qquad
y(t)=-\frac{t^6-20t^3-8}{2(t^3+2)^2}.\eqtag\label{eq.t}
$$
The cusps of~$\BB$ are certain points
$\PP_0$ over $x_0=0$,
$\PP_1$ over $x_1=1$,
and $\PP_\pm$ over $x_\pm=\Ge_\pm=(-1\pm i\sqrt3)/2$; the
corresponding values of the parameter~$t$ are
$t_0=\infty$, $t_1=1$, and $t_\pm=\Ge_\pm$. The other
points in the same fibers as the cusps correspond to the values
$t'_0=0$, $t'_1=-2$, and $t'_\pm=-2\Ge_\pm$.
The ordinate of~$\PP_0$ is $-1/2$; the ordinates of the other
three cusps are $3/2$.

The curve given by~\eqref{eq.equation} is plotted, \eg, in
Figure~\ref{fig.e6+a5} in~\S\ref{S.group}.

The curve intersects the $x$-axis at the points $x=\Ge r_\pm$,
where $r_\pm=(-1\pm\sqrt3)/2$ and $\Ge^3=1$; the corresponding
values of the parameter are $t=\Ge(1\pm\sqrt3)$. Denote the two
real intersection points by $\RR_\pm(r_\pm,0)$.

To describe the symmetries of~$\BB$, we regard them as changes of
coordinates $(x,y)$, indicating as well the corresponding change
of the parameter~$t$ in~\eqref{eq.t}. Denote by
$\Ge=(-1+i\sqrt3)/2$ a primitive cubic root of unity. Then the
group~$\AG4$ of the holomorphic automorphisms of~$\BB$ is
generated by the order~$3$ transformation
$$
(x,y)=(\Ge x',y'),\qquad t=\Ge t',\eqtag\label{eq.order3}
$$
and the order~$2$ transformation
$$
(x,y)=\Bigl(-\frac{x'-\Ge}{2\Ge^2x'+1},-\frac{3y'}{(2\Ge^2x'+1)^2}\Bigr),
\qquad
t=\frac{t'+2\Ge}{\Ge^2t'-1}.\eqtag\label{eq.order2}
$$
For the full group~$\SG4$ of Klein automorphisms of~$\BB$, one
adds
to the generating set the real structure
$\conj\:(x,y)\mapsto(\bar x,\bar y)$.

\subsection{Real structures}\label{s.conj}
In the sequel, we consider two real structures on~$\Sigma_2$ with
respect to which $\BB$ is real:
$$
\conj\:(x,y)\mapsto(\bar x,\bar y),\quad\text{and}\quad
\conj'\:(x',y')\mapsto(\bar x',\bar y'),\eqtag\label{eq.conj}
$$
where $(x',y')$ are the coordinates introduced
in~\eqref{eq.order2}
and bar stands for the complex conjugation.
The real part
$\{\Im x'=0\}$ is the circle $\ls|x+(1/2)|=3/4$, see
Figure~\ref{fig.monodromy};
it contains~$x_\pm$.
One has
$\{\Im x'=0\}\cap\{\Im x=0\}=\{\RR_+,\RR_-\}$.

Analyzing the sign of the discriminant~\eqref{eq.discrim}, one can
see that $\BB$ has three (one) $\conj$-real points over the inside
(respectively, outside) of the segment
$[x_0,x_1]$,
and it has three (one)
$\conj'$-real points over the open arc
$(x_-\,r_-\,x_+)$
(respectively, the open arc
$(x_-\,r_+\,x_+)$\,).

We are interested in real sections of~$\Sigma_2$ and their real
points. Note that sections of~$\Sigma_2$ are the parabolas of the
form
$$
y=s(x)=ax^2+bx+c.\eqtag\label{eq.section}
$$

\lemma\label{real.section}
A section~\eqref{eq.section} is real with respect to both~$\conj$
and~$\conj'$ if and only if it has the form $s(x)=-c(2x^2+2x-1)$
for some
$c\in\R$, or, alternatively, if $a$, $b$, $c$ are real and
$s(r_+)=s(r_-)=0$.
\endlemma

\proof
Under~\eqref{eq.order2}, a section~\eqref{eq.section}
transforms to $y'=a'x^{\prime2}+b'x'+c'$, where
$$
a'=\frac{a-2\Ge^2b+4\Ge c}3,\quad
b'=\frac{-2\Ge a+b+4\Ge^2c}3,\quad
c'=\frac{\Ge^2a+\Ge b+c}3.\eqtag\label{eq.abc}
$$
Equating the imaginary parts to zero, it is easy to see that all
six parameters $a$, $b$, $c$, $a'$, $b'$, $c'$ are real if and
only if $a=b=-2c$, $c\in\R$. On the other hand, $a=b=-2c$ if and
only if the section passes through~$\RR_\pm$.
\endproof

\lemma\label{real.points}
A $\conj$-real section~\eqref{eq.section} has a $\conj'$-real point
if and only if $s(r_+)s(r_-)\ge0$.
\endlemma

\proof
Assuming $a$, $b$, $c$ in~\eqref{eq.abc} real, substituting a real
value for~$x'$, and equating the imaginary part of the result to
zero, one arrives at the equation
$$
(2b+4c)x^{\prime2}-(2a+4c)x'+b-a=0.
$$
Its discriminant is $16\,s(r_+)s(r_-)$.
\endproof

\subsection{Special sections}\label{s.sections}
Let~$\BB$ be as above, and let~$\PP_0$, $\PP_1$, and~$\PP_\pm$ be
its cusps, as explained in~\ref{s.equation}.

A section~\eqref{eq.section} passes through~$\PP_0$ if and only if
$c=-1/2$; it is
tangent to~$\BB$ at~$\PP_0$ if and only if
$$
c=-1/2\quad\text{and}\quad b=0.\eqtag\label{eq.e6}
$$
Such a section passes through~$\PP_1$ if and only if
$$
(a,b,c)=(2,0,-1/2).\eqtag\label{eq.e6+a5}
$$

A section~\eqref{eq.section} passes through both~$\PP_0$
and~$\PP_1$ if and only if
$$
c=-1/2\quad\text{and}\quad a+b=2.\eqtag\label{eq.2a5}
$$
Observe that one of these sections, with $(a,b,c)=(1,1,-1/2)$, is
both $\conj$- and $\conj'$-real, see Lemma~\ref{real.section}.
The section passes through three cusps~$\PP_1$,
$\PP_\pm$ if and only if $(a,b,c)=(0,0,3/2)$.

Equating $s(x(t))=y(t)$ and $s'_t=y'_t$, one can see that
a section~\eqref{eq.section} is tangent to~$\BB$ at a
point $(x(t),y(t))\in\BB$, $t^3\ne1$,
if and only if $t\ne0$ and
$$
a=\frac{2(t^3+2)^2c+t^6+16t^3-8}{18t^2},\quad
b=-\frac{2(t^3+2)c+t^3-4}{3t}
\eqtag\label{eq.a1}
$$
or $t=0$ and $b=0$, $c=1$. Note that, if $\bar S\in\Sigma_2$
is a fixed point not over~$\PP_0$
and $t\to t_0=\infty$, then the sections as in~\eqref{eq.a1}
passing through~$\bar S$ tend to the section tangent to~$\BB$
at~$\PP_0$ (and passing through~$\bar S$).

A section as in~\eqref{eq.a1} passes through~$\PP_1$
if and only if
$$
a=\frac{2(t^3-3t-1)}{(t-1)(t+2)^2},\quad
b=\frac{6t(t+1)}{(t-1)(t+2)^2},\quad
c=-\frac{(t^3+3t^2+8)}{2(t-1)(t+2)^2}.
\eqtag\label{eq.a5+a1}
$$

Equating, in addition, $s''_t=y''_t$, one concludes that
a section~\eqref{eq.section} is inflection
tangent to~$\BB$ at a point $(x(t),y(t))\in\BB$, $t^3\ne1$,
if and only if
$$
a=\frac{t(t^3-4)}{2(t^3-1)},\quad
b=\frac{3t^2}{(t^3-1)},\quad
c=-\frac{(t^3+2)}{2(t^3-1)}.
\eqtag\label{eq.infl}
$$
Such a section cannot pass through a singular point of~$\BB$.

If a section as in~\eqref{eq.infl} is $\conj$-real,
then it has $\conj'$-real points, as one has
$$
s(r_+)s(r_-)=\frac{(t^2-2t-2)^4}{16(t^3-1)^2}\ge0,
$$
see Lemma~\ref{real.points}.
Two of the sections inflection tangent to~$\BB$
are both $\conj$- and $\conj'$-real; they
are obtained at $t=1\pm\sqrt3$
(the tangency points being~$\RR_\pm$), and their equations are
$$
y=\pm\sqrt3(2x^2+2x-1)/3.
\eqtag\label{eq.9a2}
$$
These sections are plotted in Figure~\ref{fig.9a2}
in~\S\ref{S.group}.

\subsection{Moduli and torus structures}\label{s.moduli}
The four cuspidal trigonal curve~$\BB$ has four distinct `torus
structures'; one of them is given by the decomposition
$$
2f(x,y)=8(y-3/2)^3+(6y-4x^3-5)^2,\eqtag\label{eq.torus}
$$
where $f$ is as in~\eqref{eq.equation}, and the others are
obtained from~\eqref{eq.torus} by applying a sequence of
transformations~\eqref{eq.order3} and~\eqref{eq.order2}. Formally,
a torus structure for~$\BB$ should be defined as a representation
of its equation in the form $\bp^3+e\bq^2$, where $\bp$,
$\bq$, and~$e$ are sections of $\CO_{\Sigma_2}(E+2F)$,
$\CO_{\Sigma_2}(E+3F)$, and $\CO_{\Sigma_2}(E)$, respectively,
\cf.~\cite{degt.Oka2}. Up to coefficient, $\bp$ has the form
$y-s(x)$, where $s(x)$ is a section of~$\Sigma_2$ passing through
three of the four cusps of~$\BB$. This triple of cusps determines
the torus structure. In particular, with respect to any real
structure preserving~$\BB$, two of the torus structures are real
and two are complex conjugate.
The stabilizer of each torus structure of~$\BB$ is a subgroup
$\CG3\subset\AG4$ (respectively, a dihedral subgroup
$\DG6\subset\SG4$); the stabilizer of~\eqref{eq.torus} in~$\AG4$
is~\eqref{eq.order3}.

\remark{Remark}
The four torus structures on the four-cuspidal trigonal curve~$\BB$
were also studied in Tokunaga~\cite{Tokunaga}.
\endremark

Theorem~\ref{th.reduction} implies that, in some affine
coordinates $(x,y)$ in~$\Cp2$,
each sextic~$B$ of weight eight
or nine is given by a polynomial of the form
$$
f(x, y^2+s(x)),\eqtag\label{eq.sextic}
$$
where $f$ is as in~\eqref{eq.equation} and $s(x)$ is an
appropriate section~\eqref{eq.section}, and (some of)
the torus structures
of~$B$ are obtained by substituting $y\mapsto y-s(x)$
to~\eqref{eq.torus} and its three conjugates by the automorphism
group~$\AG4$. It follows that each of the four torus
structures of a sextic of weight eight, see~\cite{degt.Oka}, is
invariant under the involution~$c$ given by
Theorem~\ref{th.reduction} and is indeed obtained
from~\eqref{eq.torus} and its conjugates.
(The twelve torus structures of a sextic of weight nine are discussed
in~\ref{s.torus.9a2} below.)

\midinsert
\table\label{tab.data}
The set of singularities of~$B$ {\it vs\.} the position of~$\LL$
\endtable
\centerline{\vbox{\offinterlineskip
\def\tm{\mathord{\times}}
\halign{\vrule height12pt depth2pt
 \quad$#$\hss\quad\vrule&
 \quad$#$\hss\quad\vrule&
 \quad\hss\smash{#}\hss\quad\vrule&
 \quad\hss\smash{$#$}\hss\quad\vrule\cr
\noalign{\hrule}
\omit\vrule height11.5pt depth3.5pt\hss Singularities\hss\vrule&
\omit\hss$\BB\cap\LL$\hss\vrule&
\omit\hss\ Condition\ \hss\vrule&
\omit\hss\ Data\ \hss\vrule\cr
\noalign{\hrule}
\bE_6\splus\bA_5\splus4\bA_2&\bA_2^*,\bA_2,\tm1&
 \eqref{eq.e6+a5}&(\PP_0,\PP_1)\cr
\bE_6\splus6\bA_2&\bA_2^*,\tm1,\tm1,\tm1&
 \eqref{eq.e6}&\{\PP_0\}\cr
2\bA_5\splus4\bA_2&\bA_2,\bA_2,\tm1,\tm1&
 \eqref{eq.2a5}&\{\PP_0,\PP_1\}\cr
\bA_5\splus6\bA_2\splus\bA_1&\bA_2,\tm2,\tm1,\tm1&
 \eqref{eq.a5+a1}&\{\PP_1\}\cr
\bA_5\splus6\bA_2&\bA_2,\tm1,\tm1,\tm1,\tm1&
 $c=-\frac12$&\{\PP_0\}\cr
8\bA_2\splus\bA_1&\tm2,\tm1,\tm1,\tm1,\tm1&
 \eqref{eq.a1}&\cr
8\bA_2&\tm1,\tm1,\tm1,\tm1,\tm1,\tm1&&\cr
\omit\vrule height2pt\hss\vrule&&&\cr
\noalign{\hrule}
9\bA_2&\tm3,\tm1,\tm1,\tm1&\eqref{eq.9a2}&\cr
\omit\vrule height2pt\hss\vrule&&&\cr
\noalign{\hrule}
\crcr}}}
\endinsert

The set of singularities of the sextic~$B$
covering~$\BB$ determines and is determined by
the singularities of the divisor $\LL+\BB$, see,
\eg,~\cite{degt.Oka3}. The relevant statements are cited in
Table~\ref{tab.data}, where the mutual position of~$\LL$ and~$\BB$
is described by listing all intersection points:
\Dashes
\roster
\dash
$\mathord{\times}n$ stands for a point of $n$-fold intersection
where $\BB$ is nonsingular,
\dash
$\bA_2$ stands for transversal intersection
of~$\LL$ and~$\BB$ at a cusp of~$\BB$, and
\dash
$\bA_2^*$ stands for a cusp of~$\BB$ where $\LL$ is tangent
to~$\BB$.
\endroster
In~\ref{s.sections}, we describe the
conditions on $(a,b,c)$ necessary for the
section~\eqref{eq.section} to be in a prescribed position with
respect to~$\BB$ (see `Condition' in the table).
If $\BB$ is fixed, in each case it
follows that the triples $(a,b,c)$ constitute a Zariski open subset in
one, four, six, or twelve irreducible families, a single family
being selected by a choice of a cusp of~$\BB$, a pair of cusps, or
an ordered pair of cusps (see `Data' in the table, where
$\{\,\cdot\,\}$ stands for a set and $(\,\cdot\,)$, for an ordered
set).

\subsection{Proof of Theorems~\ref{th.1} and~\ref{th.torus}}\label{proof.1}
Let~$B$ be a plane sextic of weight eight and with simple
singularities, let~$c$ be the involution given by
Theorem~\ref{th.reduction}, and fix an isomorphism
$\Gf\:(\Cp2(O_c)\!/c,B\!/c)\to(\Sigma_2,\BB)$. As explained in
the previous section, the equisingular moduli space of pairs
$(B,\Gf)$ consists of one to twelve connected components, a
single component being selected by a choice of a cusp of~$\BB$ or
an (ordered) pair of cusps, see Table~\ref{tab.data}. Since the
group~$\AG4$ of automorphisms of~$\BB$ acts transitively on its
cusps, pairs of cusps, and ordered pairs of cusps, ignoring~$\Gf$
results in a single deformation family, and Theorem~\ref{th.1}
follows.

For Theorem~\ref{th.torus} (in the case of simple singularities),
one should also take into consideration the torus structure
of~$\BB$, which can be identified by the cusp~$\PP$ that is
\emph{not} in the section $\{\bp=0\}$. Analyzing
Table~\ref{tab.data}, one can see that in all but one cases the
additional data still
reduce to a cusp or an (ordered) pair of cusps,
thus giving rise to a single deformation family. For example, for
the set of singularities $2\bA_5\splus4\bA_2$, the additional data
are either $(\{\PP_0,\PP_1\},\PP_0)$ or
$(\{\PP_0,\PP_1\},\PP_+)$; in the four element set of the cusps
of~$\BB$ they are equivalent to the ordered pairs $(\PP_0,\PP_1)$
and $(\PP_+,\PP_-)$, respectively. The exception is the
ordered triple $((\PP_0,\PP_1),\PP_+)$ resulting from
$(\Sigma,\Sigma_{\inj})=(\bE_6\splus\bA_5\splus4\bA_2,\bE_6\splus\bA_5\splus2\bA_2)$;
such configurations form two $\AG4$-orbits which are interchanged
by the transposition $(\PP_+,\PP_-)\in\SG4$, \ie, by the complex
conjugation.

It remains to consider pairs $(B,\{\text{\rm torus structure}\})$,
where $B$ is a sextic of weight eight with a non-simple singular
point, see~\eqref{list.J}. According to~\cite{degt.Oka2}, $B$ can
be
obtained from~$\BB$ by a birational transformation (stereographic
projection), which is determined by its blow-up center
$\OO\in\Sigma_2\sminus(\BB\cup E)$. This point either is generic
(the set of singularities $\bJ_{2,0}\splus4\bA_2$) or shares a
fiber with a cusp of~$\BB$ (the set of singularities
$\bJ_{2,3}\splus3\bA_2$). Hence, even
after an extra cusp identifying a torus
structure of~$\BB$ is added, the data selecting a connected
equisingular
deformation family still form a single $\AG4$-orbit.
\qed

\subsection{Torus structures of nine cuspidal sextics}\label{s.torus.9a2}
Pick a generic
nine cuspidal sextic $B\subset\Cp2$ and let
$C\subset\Cp2$ be the dual cubic. We identify the set
$\Sigma=\Sigma_B$ of
the cusps of~$B$ with the set of inflection points of~$C$.

A torus structure $(p,q)$ of~$B$ can be characterized by the six
point set $\Sigma_{(p,q)}\subset\Sigma$ of its inner cusps or,
equivalently, by the triple $\bSigma_{(p,q)}\subset\Sigma$ of
its outer cusps. For each cusp $P\subset\Sigma$, there are four
torus structures $(p,q)$ with $\bSigma_{(p,q)}\ni P$, and from the
discussion in~\ref{s.moduli} it follows that they are all
invariant under the involution $c_P\in\Aut(\Cp2,B)$ determined
by~$P$, see Lemma~\ref{9a2.involutions}. Hence,
each torus structure $(p,q)$ is stabilized by the
the three involutions $c_P$, $P\in\bSigma_{(p,q)}$. On the other
hand, an involution $c_P$ with $P\in\Sigma_{(p,q)}$ cannot
stabilize $(p,q)$ as any invariant subset of~$c_P$ containing~$P$
has odd cardinality. Now, using the description of $\Aut(\Cp2,B)$
given by Lemma~\ref{9a2.Aut} (and well known properties of
the inflection points on a plane cubic curve),
one can easily prove the following statements.

\proposition\label{th.3a2.torus}
If $B$ is generic, the map
$(p,q)\mapsto\{\text{\rm stabilizer}\}$ establishes a
one-to-one correspondence between the set of
twelve
torus structures of~$B$ and the set of
twelve
dihedral subgroups
$\DG6\subset\Aut(\Cp2,B)$. Each dihedral subgroup
$G\cong\DG6\subset\Aut(\Cp2,B)$
has two orbits~$\Sigma_3$ and~$\Sigma_6$ of cardinality~$3$ and~$6$,
respectively, and the inverse map sends~$G$ to the torus
structure $(p,q)$ with $\Sigma_{(p,q)}=\Sigma_6$.
\qed
\endproposition

\corollary
The action of $\Aut(\Cp2,B)$ on the set of torus structures of a
generic nine cuspidal sextic~$B$ has four orbits, each
consisting of three elements. Two torus structures
belong to the same orbit if and only if their sets of outer
singularities are disjoint.
\qed
\endcorollary

\corollary
For any nine cuspidal sextic~$B$, a triple
$\bSigma\subset\Sigma_B$ is the set of outer singularities of a
torus structure
if and only if the three points of~$\bSigma$
regarded as inflection points of the dual cubic $C\subset\Cp2$
are collinear.
\qed
\endcorollary

\corollary
The pairs $(B,\{\text{\rm torus structure}\})$, where $B$ is a
sextic of weight nine, form a connected deformation family.
\qed
\endcorollary

\section{The fundamental group}\label{S.group}

Throughout this section, $\BB\subset\Sigma_2$ is the four cuspidal
trigonal curve given by~\eqref{eq.equation}, $\LL$ is a section
of~$\Sigma_2$, and $B\subset\Cp2$ is the plane sextic obtained
as the pull-back of~$\BB$ under the double covering
$\Cp2\to\Sigma_2/E$ ramified at~$\LL$ and $E/E$.

\subsection{The braid group \rm(see~\cite{Magnus})}\label{s.Bn}
Recall that the \emph{braid group on $n$ strands}, $n\ge2$,
can be defined as the
group~$\BG{n}$ of automorphisms of the free group
$G=\<\zeta_1,\ldots,\zeta_n\>$
sending each generator to a conjugate of another generator and
leaving the product $\zeta_1\ldots\zeta_n$ fixed. We assume that the
action of~$\BG3$ on~$G$ is \emph{from the left}.
One has
$$
\BG{n}=\<\Gs_1,\ldots,\Gs_{n-1}\,|\,
 \text{$\Gs_i\Gs_{i+1}\Gs_i=\Gs_{i+1}\Gs_i\Gs_{i+1}$,
 $[\Gs_i,\Gs_j]=1$ for $\ls|i-j|>1$}\>,
$$
where
$$
\Gs_i\:(\ldots,\zeta_i,\zeta_{i+1},\ldots)\mapsto
 (\ldots,\zeta_i\zeta_{i+1}\zeta_i\1,\zeta_i,\ldots,),
\quad
i=1,\ldots,n-1,
$$
are the so called \emph{standard generators} of~$\BG{n}$.
The center of~$\BG{n}$, $n\ge3$, is
the infinite cyclic group generated by~$\Delta^2$, where
$\Delta=(\Gs_1\ldots\Gs_{n-1})(\Gs_1\ldots\Gs_{n-2})\ldots(\Gs_1)$.

%

There is a canonical homomorphism $\BG{n}\to\SG{n}$,
$\Gs_i\mapsto(i,i+1)$, where $\SG{n}$ is the symmetric group
on an $n$ element set. Its
kernel is called the \emph{pure braid group}, its elements being
\emph{pure braids}.

\subsection{Preliminary remarks}\label{s.remarks}
In this section, we briefly outline a few observations made
in~\cite{degt.Oka3} and concerning the computation of the
fundamental group of a plane sextic represented as a double of a
trigonal curve. A more formal exposition of the application of van
Kampen's method to trigonal curves can be found in~\cite{degt.e6}.

Consider the groups
$$
\pi^\Gd=\pi_1(\Sigma_2\sminus(\BB\cup E\cup\LL)),\quad
\pi^1=\pi_1(\Sigma_2\sminus(\BB\cup E)),
\quad\text{and}\quad
\pi=\pi_1(\Cp2\sminus B).
$$
To find
$\pi^\Gd$ and~$\pi^1$, we apply van Kampen's
method~\cite{vanKampen} to the vertical pencil $x=\const$ (the
ruling of~$\Sigma_2$). We choose the initial fiber~$F$ over the
point $x=r_-$, and the basis $\Ga$, $\Gd$, $\Gb$, $\Gg$ for
$\pi_F=\pi_1(F\sminus(\BB\cup E\cup\LL))$ as shown in
Figure~\ref{fig.basis}, left. (The fiber is real with respect to
both~$\conj$ and~$\conj'$, and the grey lines represent the two
real parts;
all loops are oriented in the
counterclockwise direction.) Note that $\Gd$ plays a special
r\^ole: it is a loop about $\LL\cap F$.

\midinsert
\centerline{\picture{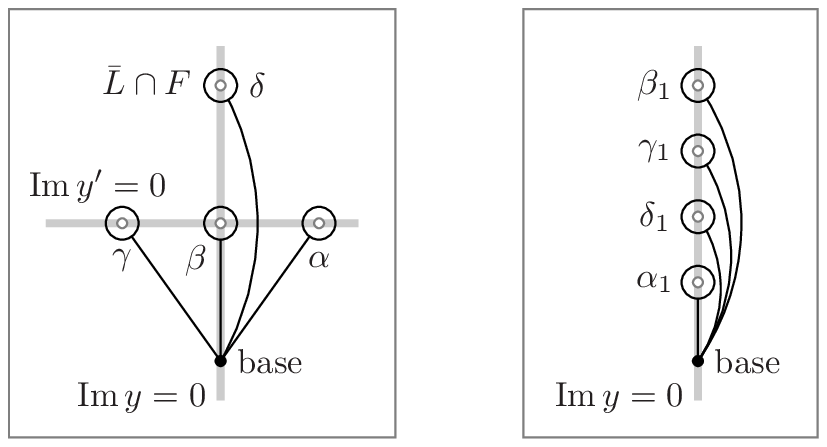}}
\figure\label{fig.basis}
The bases in the fibers
\endfigure
\endinsert

The following statement is an immediate consequence of the
double covering construction (the passage from $(\BB,\LL)$
to~$B$).

\proposition\label{d=1}
One has $\pi^1=\pi^\Gd\!/\Gd$, and
$\pi$ is the kernel of the homomorphism
$\pi^\Gd\!/\Gd^2\to\CG2$, $\Ga,\Gb,\Gg\mapsto0$, $\Gd\mapsto1$.
\qed
\endproposition

For the braid monodromy, we keep the base point in the fibers in a
$\conj$-real section $y=\const\ll0$. (Such a section is
\emph{proper} in the sense of~\cite{degt.e6}.)
The basis for
$\pi_1(\C^1\sminus\text{\{singular fibers\}})$ is chosen as shown
in Figure~\ref{fig.monodromy}
(which represents the cases resulting in sextics of weight~$8$;
the  modifications for the case of weight~$9$ are explained
in~\ref{s.w=9}). In the figure, the bold grey line
and circle represent the real parts of $\conj$ and $\conj'$,
respectively, solid (dotted) being the portions over which $\BB$
has three (respectively, one) real points. The grey dots are the
projections of the cusps of~$\BB$, and the white dot is the
projection of a point of transversal intersection of~$\BB$
and~$\LL$, see Figures~\ref{fig.e6+a5}--\ref{fig.2a5} below;
this point gives a relation $[\Gb,\Gd]=1$. For the other singular
fibers
(which all happen to be $\conj$-real and located between~$x_0$
and~$x_1$, see~\ref{s.w=8} below for details),
the loops for the monodromy
are constructed similar to~$x_1$ (see the global monodromy
computation in~\ref{s.monodromy} for more details).

\midinsert
\centerline{\picture{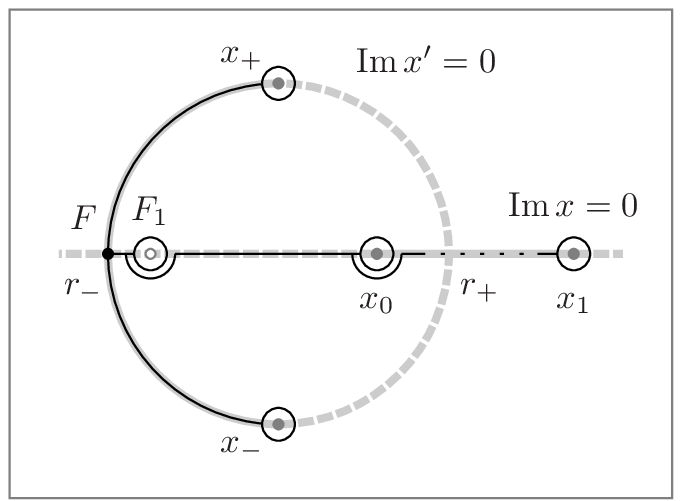}}
\figure\label{fig.monodromy}
The basis for the monodromy
\endfigure
\endinsert

Van Kampen's approach states that the group~$\pi^\Gd$ has a
presentation
$$
\pi^\Gd=\bigl<\Ga,\Gb,\Gg,\Gd\bigm|m_i=\id,\ i=1,\ldots,n,\
 (\Ga\Gd\Gb\Gg)^2=1\bigr>,\eqtag\label{eq.vanKampen}
$$
where $n$ is the number of singular fibers,
$m_i\:\pi_F\to\pi_F$ is the braid monodromy about the
singular fiber~$F_i$, and each relation $m_i=\id$
should be understood as a quadruple of relations $m_i(\Ga)=\Ga$,
$m_i(\Gb)=\Gb$, $m_i(\Gg)=\Gg$, $m_i(\Gd)=\Gd$. (For the
\emph{relation at infinity} $(\Ga\Gd\Gb\Gg)^2=1$,
see~\cite{degt.Oka3}.)

The following lemma reduces the number of singular fibers to be
considered.

\lemma\label{-1fiber}
In the presentation~\eqref{eq.vanKampen} for~$\pi^\Gd$,
\rom(any\rom) one of the braid relations $m_i=\id$ can be ignored.
\endlemma

In the sequel, we will usually ignore the monodromy about the
point $x=x_1$, as the `farthest' from~$F$.

\proof
The composition $m_1\circ\ldots\circ m_n$ of all braid monodromies
(in appropriate order) is the braid monodromy along a large circle
encompassing all singular fibers (the so called monodromy at
infinity). It is known to be $\Delta^4\in\BG4$, \ie,
the conjugation by
$(\Ga\Gd\Gb\Gg)^2$, see, \eg,~\cite{degt.e6}. Since
$(\Ga\Gd\Gb\Gg)^2=1$, the monodromy at infinity becomes the
identity and
any of~$m_i$ is a composition of the others.
\endproof

In all cases considered below, the
presentation~\eqref{eq.vanKampen}
for $\pi^\Gd$ contains a relation
$[\Gb,\Gd]=1$. Hence, the following corollary is a consequence of
Proposition~\ref{d=1} and
the standard algorithm for
presenting
a finite index subgroup, see, \eg,~\cite{MKS}.

\corollary\label{pi->pi}
If $\pi^\Gd$ is given by
$\bigl<\Ga,\Gb,\Gg,\Gd\bigm|R_j=1,\ j=1,\ldots,k\bigr>$, then
$$
\pi=\bigl<\Ga,\bGa,\Gb,\Gg,\bGg\bigm|R_j'=\bar R_j'=1,\
 j=1,\ldots,k\bigr>,
$$
where bar stands for the conjugation by~$\Gd$, $\bar w=\Gd w\Gd$,
and $R_j'$ is the relation obtained from~$R_j$, $j=1,\ldots,k$, by
eliminating~$\Gd$ \rom(using the extra relation $\Gd^2=1$\rom).
\qed
\endcorollary

\remark{Remark}
In the presentation~\eqref{eq.vanKampen} for~$\pi^\Gd$,
each relation contains an even number of copies of~$\Gd$. Hence,
after $\Gd^2=1$ is added,
each relation can be expressed in terms of
$\Ga$, $\bGa$, $\Gb$, $\Gg$, $\bGg$.
\endremark

\remark{Remark}
From Corollary~\ref{pi->pi}, the fact that the curves
used are real,
and the construction of the basis, see Figure~\ref{fig.basis},
left,
it follows that each group $G_0$, $G_1$, $G_2$, $G_3$ introduced
below has two involutive automorphisms:
$$
\gathered
\Ga\leftrightarrow\bGa,\quad
\Gb\leftrightarrow\Gb,\quad
\Gg\leftrightarrow\bGg\quad
\text{(conjugation by~$\Gd$)},\\
\Ga\leftrightarrow\Gg\1,\quad
\bGa\leftrightarrow\bGg\1,\quad
\Gb\leftrightarrow\Gb\1\quad
\text{(induced by~$\conj$)}.
\endgathered
$$
The conjugation by~$\Gd$ can as well be regarded as the
automorphism induced by the deck translation of the ramified
covering $\Cp2\to\Sigma_2/E$.
\endremark

\corollary\label{central}
If $\Gd$ is a central element of $\pi^\Gd\!/\Gd^2$, then
the map $\Ga,\bGa\mapsto\Ga$, $\Gb\mapsto\Gb$,
$\Gg,\bGg\mapsto\Gg$ establishes an isomorphism
$\pi=\pi^1$.
\qed
\endcorollary

\corollary\label{transversal}
If $\LL$ is transversal to~$\BB$, then
the map $\Ga,\bGa\mapsto\Ga$, $\Gb\mapsto\Gb$,
$\Gg,\bGg\mapsto\Gg$ establishes an isomorphism
$\pi=\pi^1$.
\endcorollary

\proof
With $\BB$ fixed, all pairs $(\BB,\LL)$ with $\LL$ transversal
to~$\BB$ are equisingular deformation equivalent. Hence, one can
take for~$\LL$ a small perturbation of $E+2F$. For such a section
it is obvious that $\Gd$ commutes with the other generators.
Indeed, in a small tubular neighborhood~$T$ of~$F$, the curve~$\BB$
is a union of three pairwise distinct almost constant sections,
whereas $\LL$ is a very `sharp' quadratic parabola
intersecting~$\BB$ transversally at six distinct points. Thus, the
space $T\sminus(\BB\cup E\cup\LL)$ is diffeomorphic to the
complement of $\C^2$ to the parabola $y=x^2$ and three distinct
horizontal lines $y=a_i\ne0$, $i=1,2,3$. It is straightforward to
see (\eg, using van Kampen's method)
that the fundamental group of the latter space has the form
$$
\bigl<\Ga_1,\Gb_1,\Gg_1,\Gd_1\bigm|
 [\Ga_1,\Gd_1]=[\Gb_1,\Gd_1]=[\Gg_1,\Gd_1]=1\bigr>,
$$
where $\Ga_1$, $\Gb_1$, $\Gg_1$, $\Gd_1$ are appropriate
generators
in~$F$
similar to those shown in Figure~\ref{fig.basis}. In particular,
the element $\Gd_1$ realized by a small loop
about
$F\cap\LL$ is in the center of the group. On the other hand,
$\Gd_1$ is conjugate to~$\Gd$, \ie, $\Gd$ is also in the center.
Then,
Corollary~\ref{central} applies to produce the desired
isomorphism.
\endproof

\remark{Remark}
Corollaries~\ref{central} and~\ref{transversal} hold for any
trigonal curve $\BB\subset\Sigma_2$,
not only the one given by~\eqref{eq.equation}.
\endremark

\subsection{The braid monodromy}\label{s.monodromy}
If the base section is chosen as indicated in~\ref{s.remarks} (a
proper section in the terminology of~\cite{degt.e6}), the braid
monodromy does indeed act on~$\pi_F$ \via\ braids.

As usual, the computation of the braid monodromy consists of two
parts. First, for a singular fiber~$F_i$, one computes the
\emph{local braid monodromy}, \ie, the braid monodromy along a
small loop about~$F_i$. For this purpose, one chooses appropriate
`standard' generators $\zeta_1,\zeta_2,\ldots$ in a nonsingular
fiber~$F_i'$
near the singular
point (similar to $\Ga,\Gb,\Gg,\ldots$ in Figure~\ref{fig.basis},
left or $\Gb_1,\Gg_1,\ldots$ in Figure~\ref{fig.basis}, right) and
uses a local standard form of the singularity to trace the points
of the curve when $F_i'$ is dragged about~$F_i$.
(We assume that $\zeta_1,\zeta_2,\ldots$ are represented by small
loops about the points of the curve that tend to the singular
point when $F_i'\to F_i$; only these
generators are involved in the braid relations resulting
from~$F_i$.)

In the present paper, we only need the singular
fibers~$F_i$ passing through the singular
points (of the curve $\BB+\LL$) of the following types:
\roster
\item
$\bA_p$, $p\ge1$: the local braid monodromy is $\Gs_1^{p+1}$
(the two points make $\frac12(p+1)$ full turns about their center of
gravity),
resulting in the only nontrivial relation
$$
\cases
(\zeta_1\zeta_2)^k=(\zeta_2\zeta_1)^k,
 &\text{if $p=2k-1$ is odd},\\
(\zeta_1\zeta_2)^k\zeta_1=(\zeta_2\zeta_1)^k\zeta_2,
 &\text{if $p=2k$ is even}.
\endcases
\eqtag\label{rel.ap}
$$
\item
$\bD_5$: the local braid monodromy is $\Gs_2^3\Gs_1\Gs_2^2\Gs_1$
(if it is $\zeta_1$ that corresponds to the nonsingular branch) or
$\Gs_1^3\Gs_2\Gs_1^2\Gs_2$ (if $\zeta_3$ corresponds to the
nonsingular branch): the two points in the singular branch make
$1.5$ full turns about their center of gravity, whereas the third
point makes one turn about the same center.
The resulting nontrivial relations are
$$
\gather
[\zeta_1,\zeta_2\zeta_3]=1,\quad
 \zeta_3\zeta_1\zeta_2\zeta_3=\zeta_1\zeta_2\zeta_3\zeta_2,
 \rlap{\quad or}
\eqtag\label{rel.d5.1}
\\
[\zeta_3,\zeta_1\zeta_2]=1,\quad
 \zeta_2\zeta_1\zeta_2\zeta_3=\zeta_1\zeta_2\zeta_3\zeta_1,\rlap{}
\eqtag\label{rel.d5.2}
\endgather
$$
respectively.
\item
$\bE_7$: the local braid monodromy is $(\Gs_1\Gs_2\Gs_1)^3$
(the two points in the singular branch make $1.5$ full turns
about the third point, which is fixed),
and the resulting nontrivial relations are
$$
[\zeta_2,\zeta_1\zeta_2\zeta_3\zeta_1]=1,\quad
 \zeta_3\zeta_1\zeta_2\zeta_3=\zeta_1\zeta_2\zeta_3\zeta_1.
\eqtag\label{rel.e7}
$$
(Here, the generator corresponding to the nonsingular branch
is~$\zeta_2$.)
\endroster

The second step is computing the \emph{global monodromy}, \ie,
relating the generators $\zeta_1,\zeta_2,\ldots$ in~$F_i'$ to the
original generators $\Ga$, $\Gb$, $\Gg$, $\Gd$ in the fixed
nonsingular fiber~$F$ and thus expressing the relations above in
terms of the original basis. For this, one chooses a path~$\xi$
connecting~$F'$ and~$F$ and drags~$F$ along~$\xi$, keeping the
base point in the base section and tracing the point of
intersection of~$F'$ and the curve. In each case considered
in the paper, both the fibers~$F$, $F'$ and the curve are real
(with respect to an appropriate real structure) and one can take
for~$\xi$ a segment of the real line circumventing the interfering
real singular fibers of the curve along small semicircles. Over
each real point of~$\xi$, all but at most two points of the curve
are real, and they can easily be traced using plots (Figures
\ref{fig.e6+a5}--\ref{fig.9a2} below); over the semicircles, the
points can be traced using local normal forms of the
singularities, similar to the computation of the local monodromy.
We leave details to the reader.

\subsection{Proof of Theorem~\ref{th.w=8}:
sextics of weight eight}\label{s.w=8}
Let $\LL$ be either the section
$y=2x^2-1/2$, see~\eqref{eq.e6+a5} and Figure~\ref{fig.e6+a5},
or its small $\conj$-real perturbation.
(In the plots, $\BB$ and~$\LL$ are shown in black and grey,
respectively.)
Such a section intersects~$\BB$ transversally at a real
point over
$x=-4/5$ (or close);
the monodromy about this point
(the fiber~$F_1$ in Figure~\ref{fig.monodromy}) gives a relation
$$
[\Gb,\Gd]=1
\eqtag\label{rel.beta}
$$
(relation~\eqref{rel.ap} with $(\zeta_1,\zeta_2)=(\Gd,\Gb)$ and
$p=1$).
Furthermore, due to Lemma~\ref{real.points},
the section has no $\conj'$-real points. Hence, over the loops
about~$x_\pm$ shown in Figure~\ref{fig.monodromy}, the point
$\LL\cap\text{\{fiber\}}$ remains above the real part, whereas the
points $\BB\cap\text{\{fiber\}}$ remain real; hence, the two
groups of points are not linked,
and one can
easily compute the monodromy. The resulting relations are
$$
\Ga\Gb\Ga=\Gb\Ga\Gb,\quad
\Gg\Gb\Gg=\Gb\Gg\Gb.\eqtag\label{rel.braid}
$$
(The $\conj'$-real part of~$\BB$ looks the same as its
$\conj$-real part, \cf. Figure~\ref{fig.e6+a5}, but the
$x'$-coordinate of~$F$ is~$r_+$ rather than~$r_-$; hence, the
relations are~\eqref{rel.ap} with $p=2$ and
$(\zeta_1,\zeta_2)=(\Ga,\Gb)$ or $(\Gb,\Gg)$, respectively.)
Finally, we have the relation at infinity
$$
(\Ga\Gd\Gb\Gg)^2=1.\eqtag\label{rel.infty}
$$
Relations~\eqref{rel.beta}--\eqref{rel.infty} hold in any
group~$\pi^\Gd$ obtained below. The corresponding relations
for~$\pi$, see Corollary~\ref{pi->pi}, are
$$
\gather
\Ga\Gb\Ga=\Gb\Ga\Gb,\quad \bGa\Gb\bGa=\Gb\bGa\Gb,
 \eqtag\label{rel.all.1}\label{rel.all.a}\\
\Gg\Gb\Gg=\Gb\Gg\Gb,\quad \bGg\Gb\bGg=\Gb\bGg\Gb,
 \eqtag\label{rel.all.2}\label{rel.all.g}\\
\Gb\Gg\Ga\Gb\bGg\bGa=1.
 \eqtag\label{rel.all.3}\label{rel.all.infty}
\endgather
$$

Letting $\Gd=1$ and adding the relation $\Ga\Gg\Ga=\Gg\Ga\Gg$
resulting from the monodromy about the cusp~$\PP_0$
(over $x=x_0$), one
obtains a presentation for the fundamental group
$G_0=\pi_1(\Sigma_2\sminus(\BB\cup E))$:
$$
G_0=\bigl<\Ga,\Gb,\Gg\bigm|
 \Ga\Gb\Ga=\Gb\Ga\Gb,\ \Gb\Gg\Gb=\Gg\Gb\Gg,\
 \Gg\Ga\Gg=\Ga\Gg\Ga,\ (\Ga\Gb\Gg)^2=1\bigr>.
\eqtag\label{eq.G0}
$$
Note that, according to~\cite{degt.Oka2}, $G_0$ is also the
fundamental group of any sextic with the set of singularities
$\bJ_{2,0}\splus4\bA_2$ or $\bJ_{2,3}\splus3\bA_2$.

\midinsert
\plot{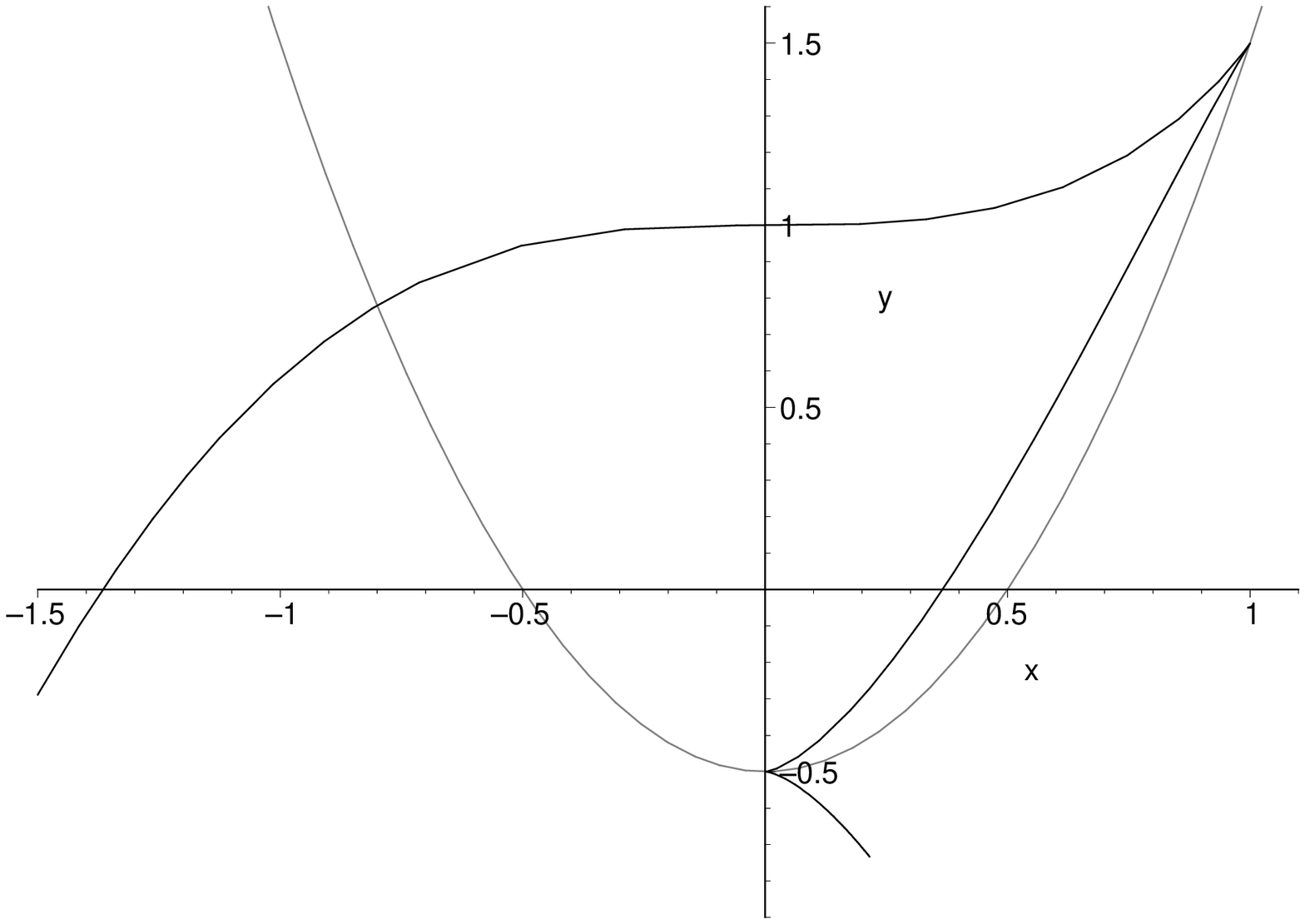}
\figure\label{fig.e6+a5}
The set of singularities $\bE_6\splus\bA_5\splus4\bA_2$
\endfigure
\endinsert

Consider the set of singularities $\bE_6\splus\bA_5\splus4\bA_2$,
see Figure~\ref{fig.e6+a5}. The additional relations for~$\pi^\Gd$
resulting from the cusp over $x=0$ are
$$
\gather
[\Gd,\Ga\Gd\Gg\Ga]=[\Gd,\Gg\Ga\Gd\Gg]=1,
\eqtag\label{rel.G2[]}
\\
 \Ga\Gd\Gg\Ga=\Gg\Ga\Gd\Gg
\eqtag\label{rel.G2=}
\endgather
$$
(relations~\eqref{rel.e7} with
$(\zeta_1,\zeta_2,\zeta_3)=(\Ga,\Gd,\Gg)$\,),
their counterpart for~$\pi$, see Corollary \ref{pi->pi},
being
$$
\bGa\Gg\Ga=\Ga\bGg\bGa=\Gg\Ga\bGg=\bGg\bGa\Gg.
\eqtag\label{rel.G2}
$$
Hence, the fundamental group is
$$
G_2=\bigl<\Ga,\bGa,\Gb,\Gg,\bGg\bigm|
 \text{\eqref{rel.all.1}--\eqref{rel.all.3}, \eqref{rel.G2}}\bigr>.
\eqtag\label{eq.G2}
$$

Note that any other curve~$B$ (or pair $\BB+\LL$)
dealt with further in this section is
a perturbation of the one just considered. Hence, the
corresponding fundamental group $\pi$ (respectively,~$\pi^\Gd$) has
relation~\eqref{rel.G2} (respectively, \eqref{rel.G2[]}
and~\eqref{rel.G2=}\,). Furthermore, similar to
Lemma~\ref{-1fiber}, in the presence of these relations one can
ignore (any) one of the singular fibers resulting from the
perturbation of the type~$\bE_7$ singular point of $\BB+\LL$ over
$x=0$.

To facilitate the further calculation, consider another real
nonsingular fiber~$F'$ between~$x=x_0$ and $x=x_1$ and the
generators $\Gb_1$, $\Gg_1$, $\Gd_1$, $\Ga_1$ for the group
$\pi_1(F'\sminus(\BB\cup E\cup\LL))$ shown in
Figure~\ref{fig.basis}, right.
Connecting~$F$ to~$F'$ by a real segment circumventing the point
$x=0$ along the semicircle $x=\epsilon e^{it}$, where $\epsilon$
is a sufficiently small positive constant and
$t\in[-\pi,0]$ (\cf. the calculation of the global monodromy
in~\ref{s.monodromy}),
and taking into account~\eqref{rel.beta},
one obtains the following expressions
$$
\Gb_1=\Ga\Gb\Ga\1,\quad
\Gg_1=\Gg,\quad
\Gd_1=\Gg\1\Gd\Gg,\quad
\Ga_1=\Gg\1\Gd\1\Ga\Gd\Gg.
\eqtag\label{eq.generators}
$$
Clearly, \eqref{eq.generators} still holds for any small
perturbation of $\BB+\LL$, as one can assume that the semicircle
above also circumvents all singular fibers resulting from the
perturbation of the type~$\bE_7$ point.
We will express the remaining relations in terms of
$\Gb_1$, $\Gg_1$, $\Gd_1$, $\Ga_1$ first, and then
use~\eqref{eq.generators} to convert them back to the original
basis.

For the set of singularities $\bE_6\splus6\bA_2$, perturb~$\LL$ to
form two new
points of transversal intersections of~$\LL$ and~$\BB$
next to and to the left from the right cusp of~$\BB$.
The new relations $[\Gd_1,\Gg_1]=[\Gd_1,\Gb_1]=1$
(relations~\eqref{rel.ap} with $p=1$ and
$(\zeta_1,\zeta_2)=(\Gg_1,\Gd_1)$ or $(\Gb_1,\Gd_1)$, respectively)
in terms of the
old generators are
$$
[\Gd,\Gg]=1\text{ (hence, $\Gd_1=\Gd$)},\quad
[\Gd,\Ga\Gb\Ga\1]=1.
$$
Since $\Ga\Gb\Ga\1=\Gb\1\Ga\Gb$ and $[\Gd,\Gb]=1$,
see~\eqref{rel.braid} and~\eqref{rel.beta}, respectively,
the last
relation implies $[\Gd,\Ga]=1$. Hence, $\Gd$ is a central element
and, due to Corollary~\ref{central},
one has $\pi=G_0$.
The same holds for the sets of singularities $\bA_5\splus6\bA_2$,
$8\bA_2\splus\bA_1$, and $8\bA_2$, which are further perturbations
of the curve just considered.

\midinsert
\plot{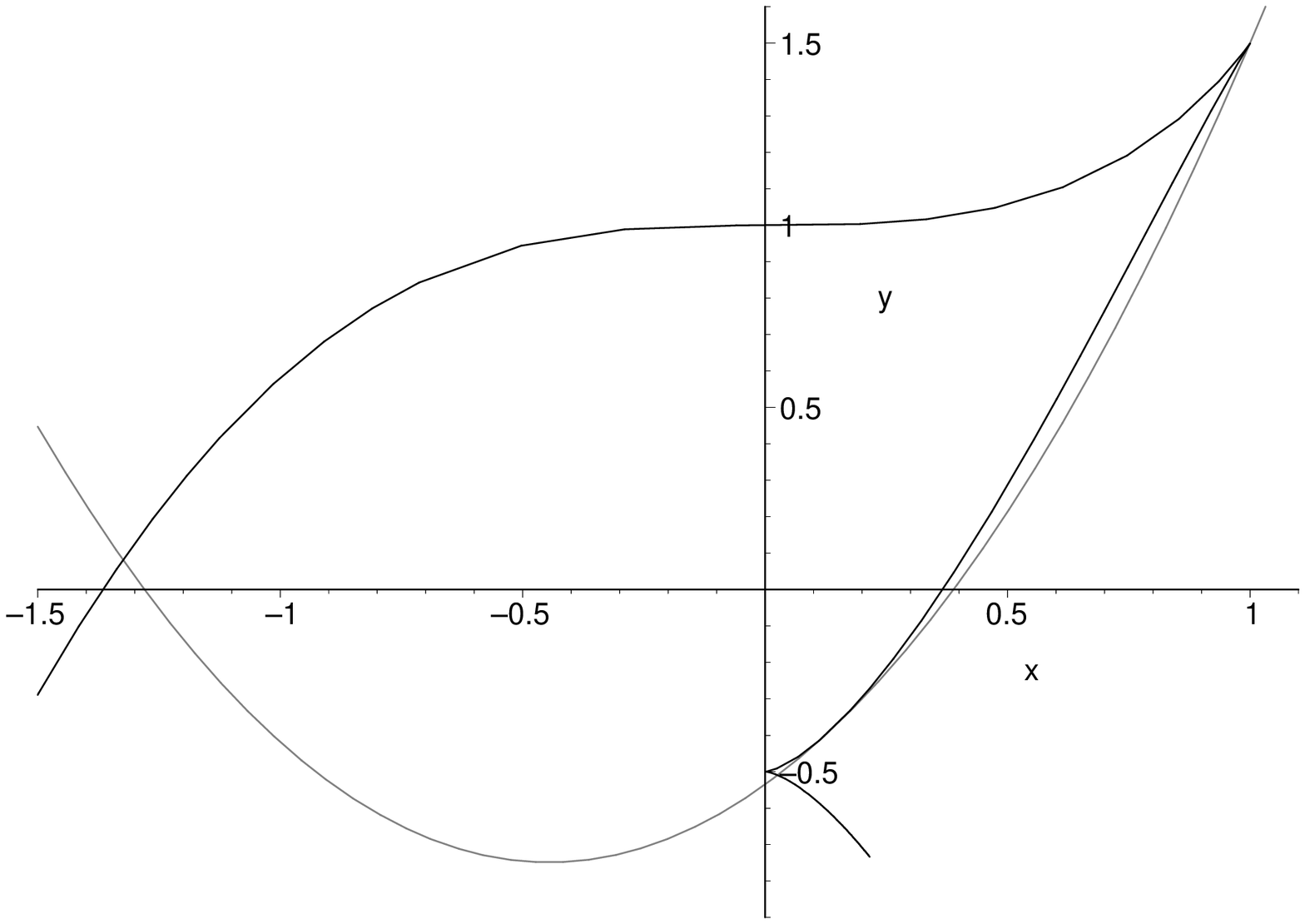}
\figure\label{fig.a5+a1}
The set of singularities $\bA_5\splus6\bA_2\splus\bA_1$
\endfigure
\endinsert

For the set of singularities $\bA_5\splus6\bA_2\splus\bA_1$,
see Figure~\ref{fig.a5+a1},
we can start with the same fiber~$F'$ and move backwards toward
the left real cusp of~$\BB$.
The extra relations are
$$
(\Gg_1\Gd_1)^2=(\Gd_1\Gg_1)^2,\quad\text{\ie,}\quad
(\Gg\Gd)^2=(\Gd\Gg)^2
$$
(from the tangency point; relation~\eqref{rel.ap} with $p=3$ and
$(\zeta_1,\zeta_2)=(\Gg_1,\Gd_1)$\,) and
$$
[(\Gg_1\Gd_1)\Gd_1(\Gg_1\Gd_1)\1,\Ga_1]=1,\quad\text{\ie,}\quad
[\Gd,(\Gd\Gg)\1\Ga(\Gd\Gg)]=1
$$
(from the point of transversal intersection; this is
relation~\eqref{rel.ap} with $p=1$ and generators
$(\zeta_1,\zeta_2)=((\Gg_1\Gd_1)\Gd_1(\Gg_1\Gd_1)\1,\Ga_1)$
obtained by circumventing the tangency point).
As explained above, in the presence of~\eqref{rel.G2[]}
and~\eqref{rel.G2=}, we can ignore the monodromy about the cusp
of~$\BB$
over $x=0$.
Since $[\Gd,(\Gd\Gg)^2]=1$, the last relation turns into
$$
[\Gd,(\Gd\Gg)\Ga(\Gd\Gg)\1]=1,\quad\text{or}\quad
[\Gd,\Gg\Ga\Gg\1]=1.
$$
Then, \eqref{rel.G2=} turns into the braid
relation $\Ga\Gg\Ga=\Gg\Ga\Gg$ and, since
$$
\Gg\Ga\Gd\Gg=(\Gg\Ga\Gg\1)\Gd\1(\Gd\Gg)^2,
$$
the commutativity relation~\eqref{rel.G2[]} becomes a
tautology. Finally, the new (in addition
to~\eqref{rel.beta}--\eqref{rel.infty}) relations for~$\pi^\Gd$
are
$$
(\Gg\Gd)^2=(\Gd\Gg)^2,\quad
\Ga\Gg\Ga=\Gg\Ga\Gg,\quad
[\Gd,\Gg\Ga\Gg\1]=1,
\eqtag\label{rel.G1.new}
$$
and their counterparts for~$\pi$, see Corollary~\ref{pi->pi}, are
$$
\Ga\Gg\Ga=\Gg\Ga\Gg,\quad
\bGa\bGg\bGa=\bGg\bGa\bGg,\quad
\Gg\Ga\bGg=\bGg\bGa\Gg,\quad
[\Gg,\bGg]=1.\eqtag\label{rel.G1}
$$
The fundamental group is
$$
G_1=\bigl<\Ga,\bGa,\Gb,\Gg,\bGg\bigm|
 \text{\eqref{rel.all.1}--\eqref{rel.all.3}, \eqref{rel.G1}}\bigr>.
\eqtag\label{eq.G1}
$$

\midinsert
\plot{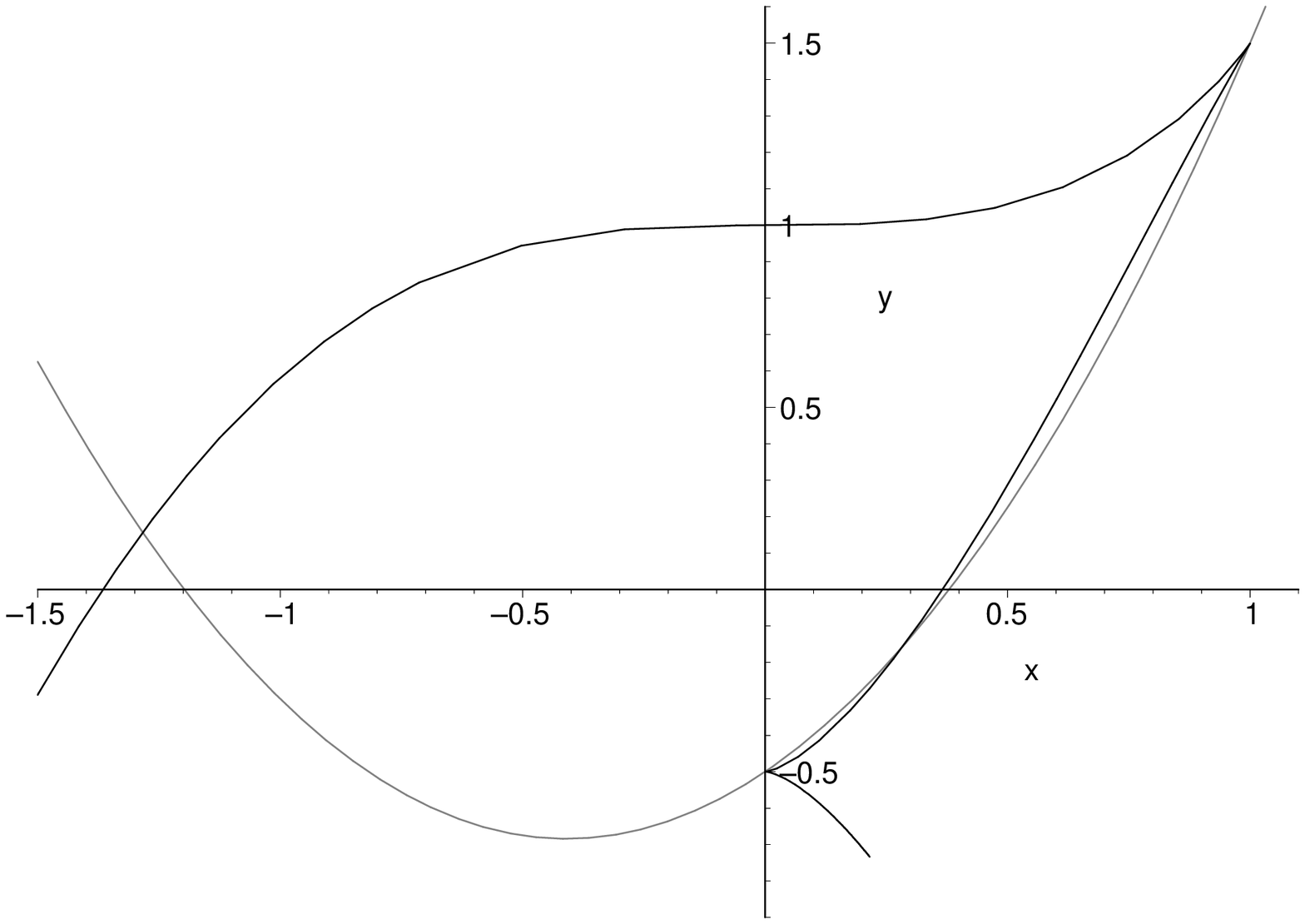}
\figure\label{fig.2a5}
The set of singularities $2\bA_5\splus4\bA_2$
\endfigure
\endinsert

Finally,
consider the set of singularities $2\bA_5\splus4\bA_2$, see
Figure~\ref{fig.2a5}. Starting from the same fiber~$F'$, one finds
that $[\Gd_1,\Gg_1]=1$ (from the point of transversal
intersection; relation~\eqref{rel.ap} with
$(\zeta_1,\zeta_2)=(\Gg_1,\Gd_1)$ and $p=1$)
and $[\Gd_1,\Gb_1\Gg_1]=[\Gd_1,\Gg_1\Ga_1]=1$
(partial relations from the two cusps; these are,
respectively, the commutativity
relations in~\eqref{rel.d5.2} and \eqref{rel.d5.1} with
$(\zeta_1,\zeta_2,\zeta_3)=(\Gb_1,\Gg_1,\Gd_1)$ and
$(\Gd_1,\Gg_1,\Ga_1)$\,). Since
$\Gb_1$, $\Gg_1$, $\Gd_1$, $\Ga_1$
also generate the group, $\Gd_1$ is a central element. Then
$\Gd=\Gd_1$, see~\eqref{eq.generators}, and since this element
is in the
center of~$\pi^\Gd$,
from Corollary~\ref{central} it
follows that
$\pi=G_0$.
\qed

\subsection{Proof of Theorem~\ref{th.w=8}:
sextics of weight nine}\label{s.w=9}
Any curve of weight nine is a nine cuspidal sextic,
see~\cite{degt.Oka}; hence, $\LL$ must be inflection tangent
to~$\BB$. We take for~$\LL$ the section
$y=\sqrt3(2x^2+2x-1)/3$, see~\eqref{eq.9a2},
which is real with respect to both~$\conj$
and~$\conj'$. This section is plotted in
Figure~\ref{fig.9a2},
where the solid grey line represents~$\LL$, and the
dotted grey line represents the same section~$\LL$ in
coordinates~$(x',y')$. (Note that $\BB$ itself looks the
same in both coordinate systems, except that the points $\RR_\pm$
trade r\^oles when the coordinates are changed.) Now, the fiber
over $x=r_-$ is singular, and we choose for~$F$ a
fiber over a $\conj$-real point between~$r_-$ and~$x_0$, \cf.
Figures~\ref{fig.monodromy} and~\ref{fig.9a2}. Then, the monodromy
about~$r_-$ gives the commutativity relation~\eqref{rel.beta}, and
hence
the monodromy about~$x_\pm$ still gives the braid
relations~\eqref{rel.braid}: due to~\eqref{rel.beta}, the
intertwining of~$\Gb$ and~$\Gd$ that occurs when $F$ moves towards
one of the two singular fibers
can be ignored.

\midinsert
\plot{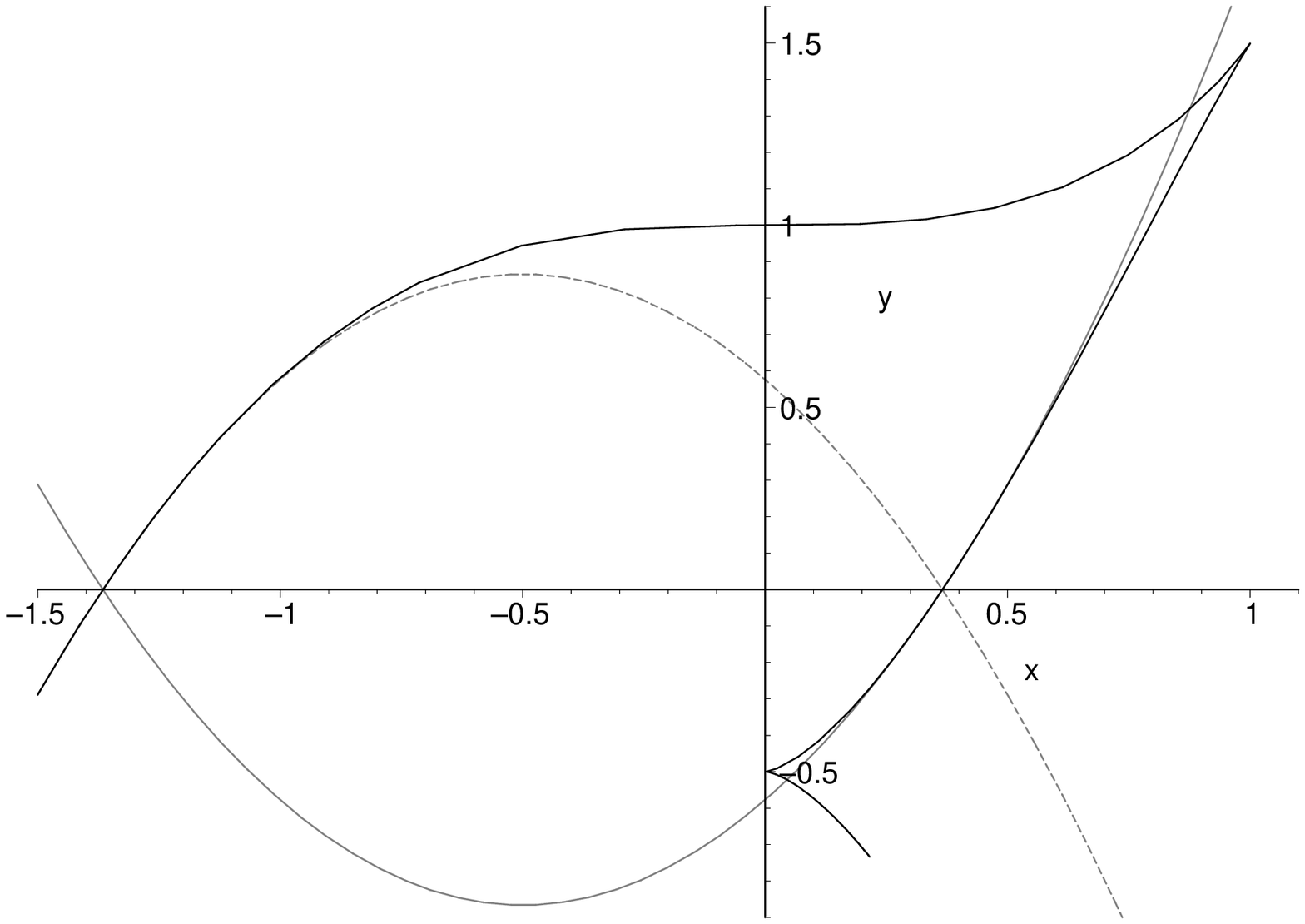}
\figure\label{fig.9a2}
The set of singularities $9\bA_2$
\endfigure
\endinsert

Consider also the generators
$\Gb_1$, $\Gg_1$, $\Ga_1$, $\Gd_1$ (\cf. Figure~\ref{fig.basis},
right) in a fiber~$F'$ over $x=\varepsilon$ for some small
$\varepsilon>0$. They can be expressed in terms of the original
generators as explained in~\ref{s.monodromy},
by circumventing the cusp over $x=0$. One has
$$
\Gb_1=\Ga\Gb\Ga\1,\quad
\Gg_1=\Gd\Gg\Gd\1,\quad
\Ga_1=(\Gd\Gg\Gd\1)\1\Ga(\Gd\Gg\Gd\1),\quad
\Gd_1=\Gd.
$$
The additional relations are
$\Gg_1\Ga_1\Gg_1=\Ga_1\Gg_1\Ga_1$ (from the cusp;
relation~\eqref{rel.ap} with $p=2$ and
$(\zeta_1,\zeta_2)=(\Gg_1,\Ga_1)$\,),
$(\Gd_1\Gg_1)^3=(\Gg_1\Gd_1)^3$ (from the inflection tangency point;
this is relation~\eqref{rel.ap} with $p=5$ and generators
$(\zeta_1,\zeta_2)=(\Gg_1,\Gd_1)$, which are
obtained by taking into account
the fact that $\Gd_1$ and $\Ga_1$ commute, see
below),
and the commutativity relations
$[\Gd_1,\Ga_1]=[(\Gg_1\Gd_1)\1\Gd_1(\Gg_1\Gd_1),\Gb_1]=1$
(from the two points of transversal intersection;
these are relations~\eqref{rel.ap} with $p=1$, where the generators
$(\zeta_1,\zeta_2)=(\Gb_1,(\Gg_1\Gd_1)\1\Gd_1(\Gg_1\Gd_1))$
at the second intersection point are obtained by circumventing the
point of inflection tangency, \cf.~\ref{s.monodromy}).
Switching back
to the original basis, we obtain
$$
\gathered
\Gd\Gg\Gd\1\Ga\Gd\Gg\Gd\1=\Ga\Gd\Gg\Gd\1\Ga,\quad
(\Gd\Gg)^3=(\Gg\Gd)^3,\\
[\Gd,(\Gd\Gg)\1\Ga(\Gd\Gg)]=
[\Gd,(\Gg\Ga)\Gb(\Gg\Ga)\1]=1.
\endgathered\eqtag\label{rel.G3.all}
$$
The corresponding relations for $\pi$, see Corollary~\ref{pi->pi},
are
$$
\gathered
\bGg\Ga\bGg=\Ga\bGg\Ga,\quad
\Gg\bGa\Gg=\bGa\Gg\bGa,\quad
\Gg\bGg\Gg=\bGg\Gg\bGg,\\
\bGg\1\Ga\bGg=\Gg\1\bGa\Gg,\quad
[\Gb,\Ga\1\Gg\1\bGg\bGa]=1,
\endgathered\eqtag\label{rel.G3}
$$
and the fundamental group is
$$
G_3=\bigl<\Ga,\bGa,\Gb,\Gg,\bGg\bigm|
 \text{\eqref{rel.all.1}--\eqref{rel.all.3}, \eqref{rel.G3}}\bigr>.\qed
\eqtag\label{eq.G3}
$$

\subsection{Proof of Theorem~\ref{th.Alexander}:
the Alexander module}\label{proof.Alexander}
Let~$B$ be a plane sextic, and let $\pi=\pi_1(\Cp2\sminus B)$ be
its fundamental group. If $B$ is irreducible, its Alexander
module~$A_B$
can be defined as follows (see Libgober~\cite{Libgober2} for
details and generalizations).

\definition\label{def.Alexander}
The \emph{Alexander module}~$A_B$ of an irreducible plane
curve~$B$ of degree~$m$ is the abelian group $\pi'\!/\pi''$
regarded as a $\Z[t]/(t^m-1)$-module \via\ the multiplication
$t[x]=[a\1xa]$, where $[x]\in A_B$ is the class of an element
$x\in\pi'$ and $a\in\pi$ is any element generating
$\pi/\pi'=\CG{m}$. The \emph{Alexander polynomial} $\Delta_B(t)$
is the order of the torsion $\C[t]/(t^m-1)$-module $A_B\otimes\C$
or, alternatively, the characteristic polynomial of the operator
$t\:A_B\otimes\C\to A_B\otimes\C$. One has
$\deg\Delta_B(t)=\rank_\Z A_B$.
\enddefinition

Thus, the Alexander module is an invariant of the fundamental
group~$\pi$. Below, we compute it for the groups $G_0$, $G_1$,
$G_2$, $G_3$ introduced in this section.

Consider the ring $\L=\Z[t]/(t^2-t+1)$.
Note that $t^3=-1$ and $t^6=1$ in~$\L$; in
particular, $t$ is invertible.
First, consider the `universal' group
$$
\tG=\bigl<\Ga,\bGa,\Gb,\Gg,\bGg\bigm|
 \text{\eqref{rel.all.1}--\eqref{rel.all.3}}\bigr>.
$$
One has $\tG\!/\tG'=\CG6$.
A standard calculation (see, \eg,~\cite{MKS})
shows that
the derived group~$\tG'$ is
generated by the elements
$$
\gathered
a_i=\Gb^i\Ga\Gb^{-i-1},\quad
c_i=\Gb^i\Gg\Gb^{-i-1},\quad i\in\Z,\\
\ba_i=\Gd a_i\Gd,\quad
\bc_i=\Gd c_i\Gd,\quad i\in\Z,\quad
\text{and}\quad b=\Gb^6,
\endgathered
$$
which are subject to the relations
$$
\gather
bs_ib\1=s_{i+6}\quad\text{for $s=a,c,\ba,\bc$ and $i\in\Z$},
 \eqtag\label{eq.G'1}\\
s_i=s_{i-1}s_{i+1}\quad\text{for $s=a,c,\ba,\bc$ and $i\in\Z$},
 \eqtag\label{eq.G'.2}\\
a_i\bc_{i+2}\ba_{i+3}c_{i+5}b=
 \ba_ic_{i+2}a_{i+3}\bc_{i+5}b=1,\quad i\in\Z.
 \eqtag\label{eq.G'.3}
\endgather
$$
Consider the abelianization $\tA=\tG'\!/\tG''$ (passing to the
additive notation) and represent~$t$ as conjugation by~$\Gb$.
One has
$s_i=t^is_0$ for $s=a,c,\ba,\bc$, $i\in\Z$, and $tb=b$.
Due to~\eqref{eq.G'.3}, $\tA$ is generated by $a_0$, $\ba_0$,
$c_0$, $\bc_0$, and from
the braid relations~\eqref{eq.G'.2} it follows that $(t^2-t+1)$
annihilates~$\tA$, \ie, $\tA$ is a $\L$-module.
Since $tb=b$, one has $b=0$ (as $t-1=t^2$ is invertible in~$\L$).
Then, in view of the fact that $t^3=-1$ in~$\L$, the last
relation~\eqref{eq.G'.3} implies
$\ba_0-a_0=t^2(\bc_0-c_0)$.
Thus, $\tA$ is a free $\L$ module with three generators,
for example, $a_0$, $c_0$, and~$\bc_0$.

Consider the group~$G_2$ corresponding to the
set of singularities
$\bE_6\splus\bA_5\splus4\bA_2$,
see~\eqref{eq.G2}. The additional relations for the
derived group~$G_2'$ are
$$
\ba_{i-1}c_ia_{i+1}=a_{i-1}\bc_i\ba_{i+1}=
 c_{i-1}a_i\bc_{i+1}=\bc_{i-1}\ba_ic_{i+1},\quad i\in\Z.
$$
Passing to the abelianization $A_2=G_2'/G_2''$ and cancelling the
term $a_0+c_1+a_2=c_0+a_1+c_2$, one obtains
$$
\ba_0-a_0=t^2(\ba_0-a_0)+t(\bc_0-c_0)=t^2(\bc_0-c_0)=
 t(\ba_0-a_0)+(\bc_0-c_0).
$$
Substituting $\ba_0-a_0=t^2(\bc_0-c_0)$, one gets
$t^2(\bc_0-c_0)=(t^3+1)(\bc_0-c_0)=0$. Thus, $\bc_0=c_0$
and $A_2$ is the free $\L$-module
generated by~$a_0$ and~$c_0$.

Any other irreducible sextic~$B$ of weight~$8$ is a perturbation of
the sextic with the set of singularities
$\bE_6\splus\bA_5\splus4\bA_2$. Hence, its Alexander
module~$A_B$ is a quotient of~$A_2$.
On the other hand,
$\rank_\Z A_B=\deg\Delta_B(t)=4$;
hence, there can be no other relations and
$A_B=A_2$.

Consider the group~$G_3$ corresponding to
the set of singularities~$9\bA_2$, see~\eqref{eq.G3}.
The extra relations for $G_3'$ are
$$
\gathered
\bc_ia_{i+1}\bc_{i+2}=a_i\bc_{i+1}a_{i+2},\quad
c_i\ba_{i+1}c_{i+2}=\ba_ic_{i+1}\ba_{i+2},\quad
\bc_ic_{i+1}\bc_{i+2}=c_i\bc_{i+1}c_{i+2},\\
\bc_i\1a_i\bc_{i+1}=c_i\1\ba_ic_{i+1},\quad
\bc_i\ba_{i+1}\ba_{i+2}\1\bc_{i+1}\1=
 c_ia_{i+1}a_{i+2}\1c_{i+1}\1,
\endgathered
$$
$i\in\Z$. Passing to the $\L$-module $A_3=G_3'/G_3''$,
the first three groups of relations
result in a tautology
$(t^2-t+1)(\,\ldots)=0$,
and the last two groups, considering that $t-1=t^2$
and $t(t-1)=-1$ in~$\L$, give the same relation
$\ba_0-a_0=t^2(\bc_0-c_0)$.
Hence, $A_3=\tA$ is the free
$\L$-module generated by $a_0$, $c_0$, and~$\bc_0$.
\qed

\section{Applications}\label{S.appl}

As above, we consider the four cuspidal
trigonal curve $\BB\subset\Sigma_2$ given by~\eqref{eq.equation}
and a section $\LL\subset\Sigma_2$ and denote by $B\subset\Cp2$
the plane sextic obtained
as the pull-back of~$\BB$ under the double covering
$\Cp2\to\Sigma_2/E$ ramified at~$\LL$ and $E/E$.

\subsection{Perturbations}\label{s.perturbations}
Denote by~$\Gamma_B$ the combined Dynkin diagram of~$B$,
\ie, the union of the Dynkin diagrams
of all singular points of~$B$ (\cf.~\ref{s.c.w=8}).

\proposition\label{subgraph}
Any induced subgraph $\Gamma'\subset\Gamma_B$
is the combined Dynkin diagram of an appropriate
small perturbation~$B'$ of~$B$. Conversely, the combined Dynkin
diagram of any perturbation of~$B$ is an induced subgraph
of~$\Gamma_B$.
\endproposition

\remark{Remark}
We do \emph{not} assert that \emph{any} sextic~$B'$ whose combined
Dynkin diagram is~$\Gamma'$ is a perturbation of~$B$. In general, this is
not true, as the homological type of~$B'$ does not need to extend
to that of~$B$.
\endremark

\remark{Remark}
Proposition~\ref{subgraph} holds for any sextic with simple
singularities; it does not need to be irreducible or
obtained by the double covering construction.
\endremark

\proof
The statement follows from the description of the moduli space of
plane sextics given in~\cite{JAG}. Let~$\Sigma_B$
and~$\Sigma'$ be the root systems spanned by the vertices of~$\Gamma_B$
and~$\Gamma'$, respectively, \cf.~\ref{s.c.w=8}. Then
$\Sigma'\oplus\<2\>\subset\Sigma_B\oplus\<2\>$ is a primitive
sublattice. Hence, it gives rise to a valid abstract homological
type, and to construct a perturbation, one shifts the
class~$\omega$
of holomorphic $2$-form of the covering $K3$-surface
from the union of the hyperplanes orthogonal
to the vertices in $\Gamma_B\sminus\Gamma'$.
Conversely, a perturbation of~$B$ means a perturbation
of~$\omega$, and we merely remove from~$\Sigma_B$ the generators
(former $(-2)$-curves) that are not orthogonal to~$\omega$.
\endproof

\corollary\label{cor.Alexander}
The Alexander polynomial of an irreducible plane sextic~$B$ of
torus type, with simple singularities, and of weight $w(B)\le7$ is
$t^2-t+1$.
\endcorollary

\proof
The Alexander polynomial can be computed using linear systems of
conics, see~\cite{poly} or~\cite{Oka.poly}, and the problem
reduces to the non-existence of a sextic of weight six (seven)
with at least two (respectively, one) singular points of weight
zero other than nodes~$\bA_1$. In view of
Proposition~\ref{subgraph}, such a sextic would perturb to a
sextic with the set of singularities $7\bA_2\splus\bA_3$. The
latter does not exist due to~\cite{JAG} or
J.-G.~Yang~\cite{Yang}: in terms of~\cite{JAG}, one would have
$\ell(\discr\tilde S)\ge5>4=\rank S^\perp$. (Here,
$S=\Sigma\oplus\Z h\subset H_2(\tX)$, \cf.~\ref{s.c.w=8},
$\tilde S$ is the
primitive hull of~$S$ in $H_2(\tX)$, and $\ell$ stands for the
minimal number of generators of a group.)
\endproof

Pick a singular point~$P$ of~$B$ and perturb it
to obtain another sextic~$B'$. To describe the impact of this
operation on the fundamental group, pick a
Milnor ball $D_P$ around~$P$ and consider a set of generators
$\Gr_1,\ldots,\Gr_k$ for the fundamental group
$\pi_1(\partial D_P\sminus B)=\pi_1(D_P\sminus B)$, the
isomorphism being induced by the inclusion
$\partial D_P\sminus B\hookrightarrow D_P\sminus B$.
After the perturbation, the inclusion homomorphism
becomes an epimorphism
$\pi_1(\partial D_P\sminus B)\twoheadrightarrow\pi_1(D_P\sminus B)$,
\ie, a number of relations $R_1,\ldots$ in the generators
$\Gr_1,\ldots,\Gr_k$ is introduced. Clearly, one has
$$
\pi_1(\Cp2\sminus B')=\pi_1(\Cp2\sminus B)/\<R_1,\ldots\>,
\eqtag\label{eq.perturbation}
$$
where the generators~$\Gr_i$ and relations $R_1,\ldots$
are regarded as elements of $\pi_1(\Cp2\sminus B)$ \via\ the
inclusion
$\partial D_P\sminus B\hookrightarrow\Cp2\sminus B$.

\subsection{Perturbations of cusps}
Let~$P$ be a cusp of~$B$ over the cusp~$\PP_+$ of~$\BB$. It can be
perturbed to either~$\bA_1$ or~$\varnothing$,
in each case the extra relation added to the
group being~$\Ga=\Gb$. (We can choose~$P$ so that $\Ga$ and~$\Gb$
generate $\pi_1(D_P\sminus B)$, and the latter group becomes
cyclic after the perturbation.)

If we perturb another cusp~$P'$ over the same
point~$\PP_+$, the group gets one more relation $\bGa=\Gb$ (and,
in view of Lemma~\ref{A2} below, this extra perturbation does not
affect the group). If we perturb another cusp~$P'$
over~$\PP_-$, the new relation is~$\Gg=\Gb$ or $\bGg=\Gb$,
depending on the choice of~$P'$; in view of Lemma~\ref{2A2} below,
the resulting groups are isomorphic (and the resulting curve is
not of torus type).

\remark{Remark}
Note that, unlike Section~\ref{S.group}, we are working with the
group of a perturbed sextic~$B'$, which
does not need to be symmetric. Hence,
we cannot automatically produce new relations applying
the conjugation by~$\Gd$.
\endremark

\lemma\label{A2}
The map $\Ga,\bGa,\Gb\mapsto\Gs_1$, $\Gg,\bGg\mapsto\Gs_2$
establishes isomorphisms
$$
G_m/\<\Ga=\Gb\>=G_m/\<\Ga=\bGa=\Gb\>=
\BG3/(\Gs_1\Gs_2\Gs_1)^2\cong\CG2*\CG3,\quad
m=0,1,2,
$$
where $G_m$ are the groups given by~\eqref{eq.G0}, \eqref{eq.G1},
and~\eqref{eq.G2}.
\endlemma

\proof
For the `smallest' group $G_0/\<\Ga=\bGa=\Gb\>$ the statement
is obvious. Hence, it suffices to consider the
`largest' group $G_2/\<\Ga=\Gb\>$.

From the relation $\Gg\Ga\bGg=\bGg\bGa\Gg$, see~\eqref{rel.G2},
it follows that
$\bGa=\bGg\1\Gg\Gb\bGg\Gg\1$. Substituting this expression and
$\Ga=\Gb$ to $\bGa\Gg\Ga=\Gg\Ga\bGg$, see~\eqref{rel.G2},
we obtain
$$
\bGg\1\Gg\cdot\underline{\Gb\bGg\Gb}=\Gg\Gb\bGg,
\quad\text{or}\quad
\bGg\1\Gg\cdot\underline{\bGg\Gb\bGg}=\Gg\Gb\bGg,
\quad\text{or}\quad
\bGg\1\Gg\bGg=\Gg.
$$
Thus, $[\Gg,\bGg]=1$ and $\bGa=\bGg\1\Gg\Gb\Gg\1\bGg$. Substitute
this expression to the braid relation $\bGa\Gb\bGa=\Gb\bGa\Gb$,
see~\eqref{rel.all.1}:
$$
\bGg\1\cdot\underline{\Gg\Gb\Gg\1}\cdot\bGg\Gb\bGg\1\cdot
 \underline{\Gg\Gb\Gg\1}\cdot\bGg=
 \Gb\bGg\1\cdot\underline{\Gg\Gb\Gg\1}\cdot\bGg\Gb.
$$
Replacing the underlined expressions using braid
relations~\eqref{rel.all.g},
one obtains
$$
\bGg\1\Gb\1\Gg\Gb\cdot\underline{\bGg\Gb\bGg\1}\cdot\Gb\1\Gg\Gb\bGg=
 \underline{\Gb\bGg\1\Gb\1}\cdot\Gg\cdot\underline{\Gb\bGg\Gb}
$$
and, once again,
$$
\bGg\1\Gb\1\Gg\Gb\Gb\1\bGg\Gb\Gb\1\Gg\Gb\bGg=
 \bGg\1\Gb\1\bGg\Gg\bGg\Gb\bGg.
$$
Cancelling and using $[\Gg,\bGg]=1$ gives $\bGg=\Gg$. Hence
also $\bGa=\Gb$, and the statement follows.
\endproof

\lemma\label{2A2}
There are isomorphisms
$$
G_m/\<\Ga=\Gb=\Gg\>=G_m/\<\Ga=\Gb=\bGg\>=\CG6,\quad
m=0,1,2,
$$
where $G_m$ are the groups given by~\eqref{eq.G0}, \eqref{eq.G1},
and~\eqref{eq.G2}.
\endlemma

\proof
Due to Lemma~\ref{A2}, all groups are quotients of
$\BG3/\<\Gs_1=\Gs_2\>=\Z$.
\endproof

\subsection{Perturbations of~$\bE_6$}
Assume that $\LL$ is tangent to~$\BB$ at~$\PP_0$ and let~$P$
be the type~$\bE_6$ singular point
over~$\PP_0$. We
consider one of the following perturbations:
$$
\bD_5,\ \bD_4,\ \bA_4\splus\bA_1,\ \bA_4,\
 \bA_3\splus\bA_1,\ \bA_3,\
 \bA_2\splus k\bA_1\ (k\le2),\
 k\bA_1\ (k\le3).
\eqtag\label{list.E6}
$$
(Thus, we exclude the perturbations $\bA_5$, $2\bA_2\splus\bA_1$,
and $2\bA_2$, which can be realized by perturbing~$\LL$ in~$\Sigma_2$
and lead to sextics of weight~$8$.)
The group $\pi_1(D_P\sminus B)$ is generated by~$\Ga$, $\bGa$, $\Gg$,
and~$\bGg$.

\lemma\label{quartic}
The space $D_P\sminus B$ is diffeomorphic to
$\Cp2\sminus(C\cup N)$, where $C\subset\Cp2$ is a plane quartic
with a type~$\bE_6$ singular point, and $N$ is a line with a
single quadruple intersection point with~$C$.
\endlemma

\proof
Since $D_P\sminus B$ is a local object, one can replace~$B$ with~$C$,
assuming that $P$ is
the type~$\bE_6$ singular point of~$C$ and that $D_P$ is a Milnor ball
about~$P$.
The affine coordinates $(x,y)$ in $\C^2=\Cp2\sminus N$ can be
chosen so that $C$ is the quartic $y^3-x^4=0$, and then
$D_P$
can be replaced with the polydisk $\Delta_1\cap\Delta_2$, where
$$
\Delta_1=\{(x,y)\in\C^2\,|\,\ls|x|<\epsilon_1\},\quad
\Delta_2=\{(x,y)\in\C^2\,|\,\ls|y|<\epsilon_2\},
$$
and
$0<\epsilon_1\ll\epsilon_2$. Then
$C\cap\Delta_1\subset\Delta_1\cap\Delta_2$, and
contracting the $y$-axis to~$\Delta_2$ results in a
diffeomorphism
$(\Delta_1\cap\Delta_2)\sminus C\cong\Delta_1\sminus C$.

Now, consider
the projection $p\:\C^2\to\C$, $(x,y)\mapsto x$, and its
restriction $p_C$ to~$C$. Each pull-back $p_C^{-1}(x)$, $x\ne0$,
consists of exactly three distinct points (as $p_C$ has no
critical points other than~$P$ and~$C$ has no branches tending to infinity
over a finite value of~$x$); hence, the fibration
$p\:\C^2\sminus C\to\C$ is locally trivial outside $x=0$, and
contracting the $x$-axis~$\C$ to~$\Delta_1$ lifts to a diffeomorphism
$\C^2\sminus C\cong\Delta_1\sminus C$.
\endproof

In view of Lemma~\ref{quartic},
the
perturbations of~$B$ inside~$D_P$ can be regarded as perturbations
of~$C$ keeping the point of quadruple intersection with~$N$,
see~\cite{quintics}, and in all cases except those excluded above,
for the perturbed quartic~$C'$ one has
$\pi_1(\Cp2\sminus(C'\cup N))=\Z$, see~\cite{groups}. Hence, the
extra relations are $\Ga=\bGa=\Gg=\bGg$. From the lemma below, it
follows that it suffices to add the relation $\Ga=\Gg$.

\lemma\label{E6}
The map $\Ga,\bGa,\Gg,\bGg\mapsto\Gs_1$, $\Gb\mapsto\Gs_2$
establishes isomorphisms
$$
G_m/\<\Ga=\Gg\>=\BG3/(\Gs_1\Gs_2\Gs_1)^2\cong\CG2*\CG3,\quad
m=0,1,2,
$$
where $G_m$ are the groups given by~\eqref{eq.G0}, \eqref{eq.G1},
and~\eqref{eq.G2}.
\endlemma

\proof
In the presence of $\Ga=\Gg$, relations~\eqref{rel.G2} are
rewritten in the form
$$
\bGa\Ga^2=\Ga\bGg\bGa=\Ga^2\bGg=\bGg\bGa\Ga.
$$
Hence, $[\Ga,\bGg\bGa]=1$, then
$\bGg=\Ga\1\bGg\bGa=\bGg\bGa\Ga\1=\bGa$, and finally
$\bGa=\bGg=\Ga=\Gg$.
\endproof

\subsection{Perturbations of~$\bA_5$}
Assume that $\LL$ passes through the cusp~$\PP_1$ of~$\BB$, and
let~$P$ be the type~$\bA_5$ singular point of~$B$ over~$\PP_1$. We
consider the following perturbations of~$P$:
$$
\gather
\bA_3\splus\bA_1,\ 3\bA_1,\eqtag\label{list.A5.2}\\
\bA_4,\ \bA_3,\ \bA_2\splus\bA_1,\ \bA_2,\ k\bA_1\ (k\le2).
\eqtag\label{list.A5.1}
\endgather
$$
(Thus, we exclude the perturbation $2\bA_2$, which can be realized
by a shift of~$\LL$.)
The group $\pi_1(D_P\sminus B)$ is generated by the
elements $\Gb_1=\Ga\Gb\Ga\1$ and~$\Gg_1=\Gg$
given by~\eqref{eq.generators}.
For the series~\eqref{list.A5.2}, the group $\pi_1(D_P\sminus B')$
is abelian of rank two; hence, the extra relation is
$[\Gb_1,\Gg_1]=1$. For the series~\eqref{list.A5.1}, one has
$\pi_1(D_P\sminus B')=\Z$ and the extra relation is
$\Gb_1=\Gg_1$.

\lemma\label{A5[]}
The map $\Ga\mapsto\Gs_1$, $\bGa\mapsto\Gs_3$, $\Gb\mapsto\Gs_2$,
$\Gg\mapsto\Gs_2\1\Gs_3\Gs_2$, $\bGg\mapsto\Gs_2\1\Gs_1\Gs_2$
establishes isomorphisms
$$
G_2/[\Gb_1,\Gg_1]=G_1/[\Gb_1,\Gg_1]=\BG4/\Gs_1^2\Gs_2\Gs_3^2\Gs_2,
$$
where $G_2$ and~$G_1$ are the groups given by~\eqref{eq.G2}
and~\eqref{eq.G1}, respectively, and $\Gb_1=\Ga\Gb\Ga\1$ and
$\Gg_1=\Gg$ are the generators given by~\eqref{eq.generators}.
\endlemma

\proof
As explained above, the relation
$[\Gb_1,\Gg_1]=1$ means that $B$ is perturbed so that the group
$\pi_1(D_P\sminus B')$
becomes abelian. On the other hand, the elements
$\Gd_1\1\Gb_1\Gd_1$ and $\Gd_1\1\Gg_1\Gd_1$ also belong to
$\pi_1(D_P\sminus B)$ and are conjugate to~$\Gg_1$ and~$\Gb_1$,
respectively. Hence, in addition, we get the relations
$$
\Gd_1\1\Gg_1\Gd_1=\Gb_1,\quad
\Gd_1\1\Gb_1\Gd_1=\Gg_1.
$$
(In view of Lemma~\ref{-1fiber}, the two additional relations should
follow from $[\Gb_1,\Gg_1]=1$ and the relations for~$G_2$, but we
will not try to deduce them algebraically.) Passing back to the
generators~$\Ga$, \dots, $\bGg$, we obtain
$$
[\Ga\Gb\Ga\1,\Gg]=1,\quad
\Ga\Gb\Ga\1=\Gg\1\bGg\Gg,\quad
\bGa\Gb\bGa\1=\bGg\1\Gg\bGg.
$$
Comparing the first two relations, we conclude that
$[\Gg,\bGg]=1$,
and then the right hand sides of the last two
relations are~$\bGg$ and~$\Gg$, respectively. Due
to~\eqref{rel.all.a}, the last two relations can be rewritten in the
form $\Gb\1\Ga\Gb=\bGg$ and $\Gb\1\bGa\Gb=\Gg$; hence,
$[\Ga,\bGa]=1$
and the relation at infinity~\eqref{rel.all.infty}
becomes $\Ga^2\Gb\bGa^2\Gb=1$.

It remains to check that all other relations follow from those
already listed and
braid relations~\eqref{rel.all.a},
so that the homomorphism in the statement is well defined.
We will do this for the smaller group $G_1/[\Gb_1,\Gg_1]$.
Braid
relations~\eqref{rel.all.g}
are the conjugation of~\eqref{rel.all.a}
by~$\Gb$, and in~\eqref{rel.G1}, after the substitution, the first two
relations turn
into~\eqref{rel.all.a}
and the third one simplifies to
$\Ga^2\Gb\bGa^2=\bGa^2\Gb\Ga^2$, which follows from
$\Ga^2\Gb\bGa^2\Gb=1$ (as both sides equal~$\Gb\1$).
\endproof

\lemma\label{A5=}
The map $\Ga,\bGa\mapsto\Gs_1$, $\Gb\mapsto\Gs_2$,
$\Gg,\bGg\mapsto\Gs_2\1\Gs_1\Gs_2$
establishes isomorphisms
$$
G_0/[\Gb_1,\Gg_1]=G_m/\<\Gb_1=\Gg_1\>=
 \BG3/(\Gs_1\Gs_2\Gs_1)^2\cong\CG2*\CG3,\quad
m=0,1,2,
$$
where $G_m$ are the groups given by~\eqref{eq.G0}, \eqref{eq.G1},
and~\eqref{eq.G2}, and $\Gb_1=\Ga\Gb\Ga\1$ and
$\Gg_1=\Gg$ are the generators given by~\eqref{eq.generators}.
\endlemma

\proof
The statement follows from Lemma~\ref{A5[]}, as we get a
relation $\Gs_3=\Gs_1$ for the braid group~$\BG4$,
either from $\Ga=\bGa$ in $G_0$, or from the relation
$\Gb_1=\Gg_1$.
\endproof

\corollary\label{no.iso}
The perturbation epimorphism $G_1\twoheadrightarrow G_0$ is not
one-to-one.
\endcorollary

\proof
If the map
$G_1\twoheadrightarrow G_0=G_1/\<\Ga\bGa\1,\Gg\bGg\1\>$
were one-to-one,
so would be its quotient
$\BG4/\Gs_1^2\Gs_2\Gs_3^2\Gs_2\twoheadrightarrow\BG3/(\Gs_1\Gs_2\Gs_1)^2$,
see Lemmas~\ref{A5[]} and~\ref{A5=}, which can be regarded as
adding an extra relation $\Gs_3=\Gs_1$. On the other hand,
$\Gs_1^2\Gs_2\Gs_3^2\Gs_2$ is a pure braid; hence, its normal
closure cannot contain
$\Gs_3\Gs_1\1$.
\endproof

\subsection{Abelian perturbations}\label{s.abelian}
Next corollary describes the result of perturbing
two or more
singular points of~$B$. (In the case of two cusps only, they
should be chosen over distinct cusps of~$\BB$.)
The statement is an immediate consequence of
Lemmas~\ref{A2}, \ref{2A2}, \ref{E6}, and~\ref{A5[]}: one should
take into account the description of the isomorphisms given by the
lemmas and observe that adding a second relation
results in the relation
$[\Gs_1,\Gs_2]=1$ in~$\BG3$, hence, in a cyclic group.

\corollary\label{2.points}
Let $G_m$, $m=0,1,2$, be one of the groups given by~\eqref{eq.G0},
\eqref{eq.G1}, or~\eqref{eq.G2}, and let $R_1$, $R_2$ be two
\emph{distinct} relations from the set $\Ga=\Gb$, $\Gg=\Gb$,
$\Ga=\Gg$, $[\Ga\Gb\Ga\1,\Gg]=1$. Then $G_m/\<R_1,R_2\>=\CG6$.
\qed
\endcorollary

\subsection{Some sextics of weight $\le7$}\label{s.w<8}
We apply the results just obtained to
produce a number of plane sextics with controlled fundamental
group. Altogether, we obtain $47$ sets of singularities of torus
type and $122$ sets of singularities not of torus type
and
not covered by Nori's theorem~\cite{Nori}.

\theorem\label{th.w=6,7}
Let $\Sigma$ be a set of singularities obtained from one of those
listed in Table~\ref{tab.torus} by removing several \rom(possibly
none\rom) nodes~$\bA_1$. Then $\Sigma$ is realized by an
irreducible plane
sextic of torus type whose fundamental group is
$\BG4/\Gs_1^2\Gs_2\Gs_3^2\Gs_2$ \rom(for the four sets of
singularities marked with a $*$ in the table\rom) or
$\BG3/(\Gs_1\Gs_2\Gs_1)^2\cong\CG2*\CG3$ \rom(for the other sets
in the table and for all nontrivial perturbations\rom).
\midinsert
\table\label{tab.torus}
Sextics of torus type
\endtable
\hbox to\hsize{\hss
\TAB\Torus
e[6] + a[5] + 3 a[2] + a[1],    [44, 7], [2, 5, 0], 18, [3, 5, false]
e[6] + a[5] + 2 a[2] + 2 a[1],  [40, 6], [3, 4, 0], 17, [4, 4, false]
e[6] + a[4] + 4 a[2],           [46, 6], [0, 5, 1], 18, [3, 5, false]
e[6] + a[3] + 4 a[2] + a[1],*   [46, 6], [2, 5, 0], 18, [3, 5, false]
e[6] + 5 a[2] + a[1],           [44, 7], [1, 6, 0], 17, [4, 6, false]
e[6] + 4 a[2] + 3 a[1],*        [42, 6], [3, 5, 0], 17, [4, 5, false]
d[5] + a[5] + 4 a[2],           [47, 6], [2, 5, 0], 18, [3, 5, false]
d[5] + 6 a[2],                  [47, 6], [1, 6, 0], 17, [4, 6, false]
d[4] + a[5] + 4 a[2],           [45, 6], [3, 5, 0], 17, [4, 5, false]
d[4] + 6 a[2],                  [45, 6], [2, 6, 0], 16, [5, 6, false]\ENDTAB
\hss
\TAB\Torus
2 a[5] + 3 a[2] + a[1],         [44, 7], [3, 5, 0], 17, [4, 5, false]
2 a[5] + 2 a[2] + 2 a[1],       [40, 6], [4, 4, 0], 16, [5, 4, false]
a[5] + a[4] + 4 a[2] + a[1],    [48, 6], [2, 5, 1], 18, [3, 5, false]
a[5] + a[3] + 4 a[2] + a[1],    [46, 6], [3, 5, 0], 17, [4, 5, false]
a[5] + 5 a[2] + 2 a[1],         [46, 7], [3, 6, 0], 17, [4, 6, false]
a[5] + 4 a[2] + 3 a[1],         [42, 6], [4, 5, 0], 16, [5, 5, false]
a[4] + 6 a[2] + a[1],           [48, 6], [1, 6, 1], 17, [4, 6, false]
a[3] + 6 a[2] + 2 a[1],*        [48, 6], [3, 6, 0], 17, [4, 6, false]
7 a[2] + 2 a[1],                [46, 7], [2, 7, 0], 16, [5, 7, false]
6 a[2] + 4 a[1],*               [44, 6], [4, 6, 0], 16, [5, 6, false]\ENDTAB\hss}
\endinsert
\endtheorem

\proof
We start with one of the sets of singularities listed
in~\eqref{list.w=8}, realize it by a sextic~$B$ of weight
eight, and use Proposition~\ref{subgraph} to
perturb~$B$ to a new sextic~$B'$.
The fundamental group of~$B'$ is
controlled using Lemmas~\ref{A2},
\ref{E6}, \ref{A5[]}, and~\ref{A5=}. According to the lemmas, one
can perturb up to two cusps of~$B$ (both over the same cusp
of~$\BB$) or one type~$\bE_6$ or~$\bA_5$ singular point of~$B$,
so that the new group $\pi_1(\Cp2\sminus B')$ would factor to
$\CG2*\CG3$. (In the case of a type~$\bE_6$
or~$\bA_5$ singular point, the allowed perturbations are those listed
in~\eqref{list.E6} or~\eqref{list.A5.2}, \eqref{list.A5.1},
respectively.)
Since the latter group has nontrivial Alexander
polynomial $t^2-t+1$, the new sextic~$B'$ is of torus type,
see~\cite{degt.Oka}.
\endproof

\remark{Remark}
Theorem~\ref{th.w=6,7} covers five tame sextics:
$$
\bE_6\splus\bA_5\splus2\bA_2,\quad
\bE_6\splus4\bA_2,\quad
2\bA_5\splus2\bA_2,\quad
\bA_5\splus4\bA_2,\quad
6\bA_2.
$$
The groups of these curves were first found in~\cite{OkaPho}.
\endremark

\corollary
The Alexander module of each curve mentioned in
Theorem \ref{th.w=6,7} is a free module on one generator over the
ring $\L=\Z[t]/(t^2-t+1)$.
\endcorollary

\proof
The Alexander module of each braid group~$\BG{n}$, $n\ge3$, is
known to be a free $\L$-module generated, \eg, by $\Gs_2\Gs_1\1$.
Since the Alexander polynomial of each curve is $t^2-t+1$, there
can be no further relations.
\endproof

\theorem\label{th.nontorus}
Let $\Sigma$ be a set of singularities obtained from one of those
listed in Table~\ref{tab.nontorus} by several \rom(possibly
none\rom) perturbations $\bA_2\to\bA_1,\varnothing$ or
$\bA_1\to\varnothing$.
Then $\Sigma$ is realized by an irreducible plane
sextic, not of torus type, whose fundamental group is~$\CG6$.
\midinsert
\table\label{tab.nontorus}
Sextics with abelian fundamental group
\endtable
\hbox to\hsize{\hss
\TAB\NonTorus
e[6] + a[5] + 2 a[2] + 2 a[1],  [40, 6], [3, 4, 0], 17, [4, 4, false]
e[6] + a[4] + 3 a[2] + a[1],    [42, 5], [1, 4, 1], 17, [4, 4, false]
e[6] + a[3] + 3 a[2] + 2 a[1],  [42, 5], [3, 4, 0], 17, [4, 4, false]
e[6] + 4 a[2] + 2 a[1],         [40, 6], [2, 5, 0], 16, [5, 5, false]
e[6] + 3 a[2] + 4 a[1],         [38, 5], [4, 4, 0], 16, [5, 4, false]
d[5] + a[5] + 3 a[2] + a[1],    [43, 5], [3, 4, 0], 17, [4, 4, false]
d[5] + a[4] + 4 a[2],           [45, 4], [1, 4, 1], 17, [4, 4, false]
d[5] + a[3] + 4 a[2] + a[1],    [45, 4], [3, 4, 0], 17, [4, 4, false]
d[5] + 5 a[2] + a[1],           [43, 5], [2, 5, 0], 16, [5, 5, false]
d[5] + 4 a[2] + 3 a[1],         [41, 4], [4, 4, 0], 16, [5, 4, false]
d[4] + a[5] + 3 a[2] + a[1],    [41, 5], [4, 4, 0], 16, [5, 4, false]
d[4] + a[4] + 4 a[2],           [43, 4], [2, 4, 1], 16, [5, 4, false]
d[4] + a[3] + 4 a[2] + a[1],    [43, 4], [4, 4, 0], 16, [5, 4, false]\ENDTAB
\hss
\TAB\NonTorus
d[4] + 5 a[2] + a[1],           [41, 5], [3, 5, 0], 15, [6, 5, false]
d[4] + 4 a[2] + 3 a[1],         [39, 4], [5, 4, 0], 15, [6, 5, false]
2 a[5] + 2 a[2] + 2 a[1],       [40, 6], [4, 4, 0], 16, [5, 4, false]
a[5] + a[4] + 3 a[2] + 2 a[1],  [44, 5], [3, 4, 1], 17, [4, 4, false]
a[5] + a[3] + 3 a[2] + 2 a[1],  [42, 5], [4, 4, 0], 16, [5, 4, false]
a[5] + 4 a[2] + 3 a[1],         [42, 6], [4, 5, 0], 16, [5, 5, false]
2 a[4] + 4 a[2] + a[1],         [46, 4], [1, 4, 2], 17, [4, 4, false]
a[4] + a[3] + 4 a[2] + 2 a[1],  [46, 4], [3, 4, 1], 17, [4, 4, false]
a[4] + 5 a[2] + 2 a[1],         [44, 5], [2, 5, 1], 16, [5, 5, false]
a[4] + 4 a[2] + 4 a[1],         [42, 4], [4, 4, 1], 16, [5, 4, false]
2 a[3] + 4 a[2] + 2 a[1],       [44, 4], [4, 4, 0], 16, [5, 4, false]
a[3] + 5 a[2] + 3 a[1],         [44, 5], [4, 5, 0], 16, [5, 5, false]
6 a[2] + 3 a[1],                [42, 6], [3, 6, 0], 15, [6, 6, false]
5 a[2] + 5 a[1],                [40, 5], [5, 5, 0], 15, [6, 5, false]\ENDTAB\hss}
\endinsert
\endtheorem

\proof
As in the previous proof, we start with a sextic~$B$ of weight
eight (one of the sets of singularities listed
in~\eqref{list.w=8}\,)
and use Proposition~\ref{subgraph} to perturb two or more singular
points of~$B$. (If only two cusps are perturbed, they should be
over distinct cusps of~$\BB$; the allowed perturbations
of a type~$\bE_6$
or~$\bA_5$ singular point are those listed
in~\eqref{list.E6} or~\eqref{list.A5.2}, \eqref{list.A5.1},
respectively.)
Then, due to Corollary~\ref{2.points},
the new group is abelian.
\endproof

\remark{Remark}
Since we only construct examples in the proofs,
Theorems~\ref{th.w=6,7} and~\ref{th.nontorus} are stated in the
form of existence. Proving the uniqueness of the equisingular
deformation family realizing each of the $169$ sets of
singularities covered by the theorem would require a lot of
tedious calculations. (Certainly, when speaking about uniqueness,
one should distinguish sextics of torus type and those not of
torus type; in some cases of weight~$6$,
see item~\therosteritem2 below,
they may share the same sets of singularities.) We mention two
series where the uniqueness is known.
\roster
\item\local{JAG}
For $19$ sets of
singularities in Theorem~\ref{th.w=6,7} and for $60$ sets of
singularities in Theorem~\ref{th.nontorus}, the uniqueness of an
equisingular deformation family follows directly from
Theorem~5.2.1 in~\cite{JAG}.
\item\local{Aysegul}
According to A.~\"Ozg\"uner~\cite{Aysegul}, the uniqueness
takes place for any set of singularities of weight~$6$ with all
points of weight zero of type~$\bA_1$. With few exceptions,
these sets of
singularities are realized by two deformation families: one of
torus type and one not.
\endroster
\endremark

Comparing the two lists and using~\cite{Aysegul},
we obtain the following corollary
related to some of the so called \emph{classical Zariski pairs},
\ie, pairs of sextics
that share the same set of singularities but differ by the
Alexander polynomial, see~\cite{poly}.

\corollary\label{cor.Zariski}
Each of the following $17$ sets of singularities
$$
\let\1\quad\def\2#1{\gathered#1\endgathered}\def\3#1{\quad(k\le#1)}\setcat\2{
e[6] + a[5] + 2 a[2] + k a[1],\1
e[6] + 4 a[2] + k a[1],\1
2 a[5] + 2 a[2] + k a[1]\32,\\
a[5] + 4 a[2] + k a[1],\1
6 a[2] + k a[1]\33}
$$
is realized by exactly two equisingular deformation families of
irreducible plane sextics\rom. In each case,
the families differ by the fundamental group of the curves, which
is either $\CG2*\CG3$ \rom(torus type\rom) or~$\CG6$.
\qed
\endcorollary

\refstyle{A}

\enddocument